\documentclass[a4paper, 11pt]{amsart}
\usepackage{amsmath,amsthm,amssymb,amscd,amsfonts}
\usepackage{tikz-cd}
\usepackage{enumitem}

\usepackage[backref=page]{hyperref}

\newtheorem{fact}{Theorem}

\newtheorem{theorem}{Theorem}[section]
\newtheorem{proposition}[theorem]{Proposition}
\newtheorem{lemma}[theorem]{Lemma}
\newtheorem{corollary}[theorem]{Corollary}
\theoremstyle{remark}
\newtheorem{remark}{Remark}
\newtheorem{example}{Example}

\def\Z{\mathbb Z}
\def\Q{\mathbb Q}
\def\R{\mathbb R}
\def\C{\mathbb C}

\def\pt{{\rm pt}}

\begin{document}

\title{On bordism classes of symmetric products of surfaces}

\author{G. S. Chernykh}
\address{Steklov Mathematical Institute of the Russian Academy of Sciences, Moscow, Russia;\newline
National Research University Higher School of Economics, Russian Federation}
\email{aaa057721@gmail.com}

\author{D. V. Gugnin}
\address{Faculty of Mechanics and Mathematics, Lomonosov Moscow
State University, Russia}
\email{dmitry-gugnin@yandex.ru}

\thanks{The study has been funded within the framework of the HSE University Basic Research Program  (HSE-BR-2025-023). The first author is a winner of the Theoretical Physics and Mathematics Advancement Foundation “BASIS” program for postdocs.}

\begin{abstract}
We explicitly determine the complex bordism class of the $n$-th symmetric product of a complex curve of genus $g$. In particular, we provide an explicit formula expressing this class as an integral linear combination of bordism classes of products of complex projective spaces.

Similarly, we determine the unoriented bordism class of the $n$-th symmetric product of a non-orientable surface of genus $g$ (i.e., $\#_{i=1}^{g+1}\R P^2$). It turns out that these classes are also expressed as a $\Z_2$-linear combination of bordism classes of products of even-dimensional real projective spaces.

We also give an affirmative answer to a non-orientable version of the Blagojevi\'c--Gruji\'c--\v{Z}ivaljevi\'c conjecture stating that the topological type of any symmetric product of a noncompact surface of finite type determines the topological type of the surface itself.
\end{abstract}

\maketitle

\section{Introduction}

For a Hausdorff space $X$, its \emph{$n$-th symmetric product} $SP^n(X)$ is the quotient $$(\underbrace{X\times \cdots \times X}_{n \mbox{ {\scriptsize times}}})/\Sigma_n$$ of the $n$-fold Cartesian product by the action of the permutation group. One can easily check that $SP^n(X)$ is also Hausdorff. The elements of $SP^n(X)$ are represented by (unordered) multisets of points of $X$. We will denote them by $[x_1, \ldots, x_n]$. Moreover, any homotopy $f_t\colon X \to Y$ induces a homotopy $SP^n(f_t)\colon SP^n(X)\to SP^n(Y)$. So the homotopy type of $SP^n(X)$ depends only on the homotopy type of $X$ (and the number $n$).

The theory of symmetric products is a classical subject of study in geometry and topology (see, for instance, the nice survey \cite{bgz}). From the topological point of view, since the mid-20th century, important results concerning the homotopy and (co)homology properties of symmetric products have been obtained; see, for example, the classical references \cite{nak, dt, dp, milg}, and also the more recent \cite{br} (where the integral cohomology ring of the symmetric squares of complex and quaternionic projective spaces is calculated) and \cite{j25} (where the topological complexity of symmetric products of surfaces is computed). Also, see \cite{gug4}, where the integral cohomology ring modulo torsion of symmetric products $SP^n(X)$ is described (in terms of the integral cohomology ring modulo torsion of the original $X$ itself). More recently, symmetric products are playing an important role in the theory of $n$-valued multiplications on spaces ($n$-valued topological groups, $nH$-spaces, etc.), see \cite{buch06, gug1, gug2, gug3}.

There is a natural homeomorphism $SP^n(\C)\cong \C^n$ sending the multiset $[z_1,\ldots, z_n]$ to the ordered tuple of coefficients of the polynomial $\prod_{i=1}^n(t - z_i)$. In particular, $SP^n(\C)$ acquires a natural complex structure.

In fact, it is well known that for any Riemann surface $C$, from these local charts one can obtain a natural complex $n$-dimensional manifold structure on $SP^n(C)$ (essentially because $\mathcal O_{\C^n}^{\Sigma_n} \cong \mathcal O_{\C^n}$ for the sheaf of holomorphic functions). For example, $SP^n(\C P^1)\cong\C P^n$ (via a similar biholomorphism sending the multiset $[(z_1:w_1),\ldots, (z_n:w_n)]$ to the ordered tuple of coefficients of the polynomial $\prod_{i=1}^n(z_i t - w_i)$).

More generally, for any smooth real surface $N$, we similarly have that $SP^n(N)$ is a $2n$-dimensional smooth manifold. For example, it is known that $SP^n(\R P^2)\cong \R P^{2n}$ (see \cite{dulu}, or \cite{arn}). For the convenience of the reader, we add details on the existence and uniqueness of such smooth structures in the Appendix (where we also consider the question of which subgroups $G\subset \Sigma_n$ make the quotient $M^{\times n}/G$ a manifold).

On the other hand, for a compact Riemann surface $C$, its symmetric product $SP^n(C)$ is a projective variety whose algebraic geometry has been studied in a series of papers; for example, see the classical ones \cite{mat1, mat2, mat3, mcd, sch}.

In this paper, we primarily consider symmetric products of surfaces from the bordism point of view. Namely, we generalize Macdonald's calculations \cite{mcd} of the signature and the $\chi_y$-genus of symmetric products of complex curves and give explicit formulae for complex (or unoriented) bordism classes of symmetric products of closed surfaces. It turns out that they are expressed as linear combinations of products of projective spaces, i.e., they lie in specific subrings of the corresponding bordism rings generated by projective spaces.

Since a symmetric product of a surface is an irreducible projective variety, our results can also be considered in the context of a general Hirzebruch problem: which bordism classes contain irreducible projective representatives?

On the other hand, from a topological point of view, symmetric products of Riemann surfaces were studied, for example, in \cite{bgz1, bgz2}, where the authors conjectured that the topological type of any symmetric product of a Riemann surface determines the topological type of the surface itself. In \cite{gug4} this conjecture was proved for oriented surfaces. In this paper we also extend this result to non-orientable surfaces.

The authors are grateful to Victor Buchstaber for posing the problem and for fruitful discussions.

\section{Main results}

For a (closed) complex or real manifold $M$, we will denote by $[M]$ its bordism class in the complex bordism ring $\varOmega^U_*$ or in the unoriented bordism ring $\varOmega^O_*$, respectively (the details of these bordism theories can be found in \cite[Chapters VI, VII]{sto}). Whether we are in the complex or the real case will be clear from the context.

Recall that by Thom's theorem \cite{thom2} the ring $\varOmega^O_*$ is a polynomial algebra over $\Z_2$ with one generator in each dimension not of the form $2^k-1$. Moreover, the even-dimensional real projective spaces can be chosen as the polynomial generators in even dimensions (odd-dimensional projective spaces are null-bordant).

   Also, by the Milnor--Novikov theorem \cite{mil}, \cite{nov} the ring $\varOmega^U_*$ is polynomial over $\Z$ with generators in each even dimension, and complex projective spaces can be chosen as polynomial generators of $\varOmega^U_*\otimes\Q$.  However, over the integers $\C P^n$ can be chosen as a polynomial generator if and only if $n+1$ is prime.

For a smooth closed (not necessarily connected) complex curve $C$, consider the following formal power series
$$SP(C,t)=1+\sum_{n=1}^{\infty} [SP^n(C)]t^n \in \varOmega^U_*[[t]].$$

For example, $SP(\C P^1, t) = 1 + \sum [\C P^n]t^n =: \C P(t)$.

Similarly, for a smooth closed (not necessarily connected) real surface $N$, we consider the formal power series
$$SP(N,t)=1+\sum_{n=1}^{\infty} [SP^n(N)]t^n \in \varOmega^O_*[[t]].$$ We have $SP(\R P^2, t) = 1 + \sum [\R P^{2n}]t^n =: \R P(t)$.

Note that $C$ and $N$ above are not assumed to be connected, and we have the formula \begin{equation}\label{sum}
    SP(A\sqcup B,t)=SP(A,t)\cdot SP(B,t)
\end{equation}

Notice also that, since the series $\C P(t)$ and $\R P(t)$ have constant term $1$, there exist power series $\frac{1}{\C P(t)}$ and $\frac{1}{\R P(t)}$.

Our main result is the following.

\begin{theorem}\label{th1}
    \begin{enumerate}[label=(\alph*)]
    \item Let $C_g$ be a complex curve of genus $g$. Then $$SP(C_g,t)=\Bigl(\frac{1}{\C P(t)}\Bigr)^{g-1}.$$

    \item Let $N_g$ be the connected sum $\underbrace{\R P^2\#\cdots\#\R P^2}_{g+1}$. Then $$SP(N_g,t)=\Bigl(\frac{1}{\R P(t)}\Bigr)^{g-1}.$$
    \end{enumerate}
\end{theorem}

Comparing coefficients of $t^n$ in the expressions above, we obtain the following explicit formulae, expressing bordism classes (complex or unoriented) of symmetric products of surfaces as integral linear combinations of bordism classes of products of projective spaces:

\begin{equation}\label{fu}
    [SP^n(C_g)]=\sum\limits_{|\omega|=n} (-1)^{\ell(\omega)} k(\omega) \binom{g+\ell(\omega)-2}{\ell(\omega)} [\C P^\omega],
\end{equation}

\begin{equation}\label{fo}
    [SP^n(N_g)]=\sum\limits_{|\omega|=n} (-1)^{\ell(\omega)} k(\omega) \binom{g+\ell(\omega)-2}{\ell(\omega)} [\R P^{2\omega}],
\end{equation}
where $\omega=(i_1 \ge \ldots \ge i_s\ge 1)$ is an integer partition of a number $|\omega|=i_1+\ldots + i_s$, $\ell(\omega)=s$ is its length, $\C P^\omega = \C P^{i_1} \times \cdots \times \C P^{i_s}$, $\R P^{2\omega} = \R P^{2i_1} \times \cdots \times \R P^{2i_s}$, $\binom{n}{k}$ is a binomial coefficient (zero if $n<k$) and $$k(\omega)=\#\{\text{distinct ordered tuples }(i_{\sigma(1)}, \ldots, i_{\sigma(s)}) \text{ obtained by permuting the entries of }\omega\}.$$ Denote by $\ell_1(\omega)=\#\{1\le r \le s \mid i_r = 1\}$, $\ell_2(\omega)=\#\{1\le r \le s \mid i_r = 2\}$, etc. Then $\ell(\omega)=\ell_1(\omega)+\ell_2(\omega)+\ldots$ and  $k(\omega)=\frac{\ell(\omega)!}{\ell_1(\omega)! \ell_2(\omega)! \ldots}$.

\begin{remark}
    Of course, in formula \eqref{fo} we consider coefficients modulo $2$; in particular, signs do not matter.
\end{remark}

\begin{proof}[Deduction of formulae \eqref{fu} and \eqref{fo} from Theorem \ref{th1}]
Let $1+A(x)=1+\sum_{n\ge 1} a_n x^n$.

    We have the Taylor expansion $$(1+x)^{-n}=1+\sum_{k\ge 1} \frac{(-n)(-n-1)\cdots(-n-k+1)}{k!} x^k=1+\sum_{k\ge 1} (-1)^k\binom{n+k-1}{k}x^k$$

    So, $${\rm coeff}_n \frac{1}{(1+A(x))^{g-1}}=\sum_{k\ge 1} (-1)^k\binom{n+k-1}{k}{\rm coeff}_n\, A(x)^k$$

    Clearly, $${\rm coeff}_n\, A(x)^k=\sum_{\substack{i_1+\cdots+i_k=n\\ i_j\ge 1}} a_{i_1}\cdot\ldots \cdot a_{i_k}$$

    Now, grouping like terms (ordering $i_j$), we obtain $${\rm coeff}_n (A(x)^k)=\sum_{\substack{|\omega|=n\\ \ell(\omega)=k}} k(\omega)a_\omega$$
\end{proof}

 So, Theorem \ref{th1} says that symmetric products of surfaces lie in specific subrings of the bordism rings generated by projective spaces (i.e., generated ``by half'' of the generators in the unoriented case, and by rational, but not integral, generators in the complex case).

 \begin{example}
     We give some explicit examples of the formulae above.

     $$SP^2(C_g)=-(g-1)[\C P^2]+\frac{g(g-1)}{2}[\C P^1]^2$$

     $$SP^3(C_g)=-(g-1)[\C P^3]+g(g-1)[\C P^2][\C P^1]-\frac{(g+1)g(g-1)}{6}[\C P^1]^3$$

     \begin{multline*}
        SP^4(C_g)= -(g-1)[\C P^4]
+\frac{g(g-1)}{2}[\mathbb{C}P^2][\mathbb{C}P^2]
+g(g-1)[\mathbb{C}P^3][\mathbb{C}P^1]-\\
-\frac{(g+1)g(g-1)}{2}[\C P^2][\mathbb{C}P^1]^2+
\frac{(g+2)(g+1)g(g-1)}{24}[\mathbb{C}P^1]^4
     \end{multline*}

     \begin{multline*}
         SP^5(C_g)=-(g-1)[\C P^5]
+ g(g-1)[\C P^3][\C P^2]
+ g(g-1)[\C P^4][\C P^1]-\\
- \frac{(g+1)g(g-1)}{2}[\C P^2]^2[\C P^1]
- \frac{(g+1)g(g-1)}{2}[\C P^3][\C P^1]^2+\\
+ \frac{(g+2)(g+1)g(g-1)}{6}[\C P^2][\C P^1]^3
- \frac{(g+3)(g+2)(g+1)g(g-1)}{120}[\C P^1]^5
\end{multline*}

 \end{example}

\medskip
Before proving Theorem \ref{th1}, let us formulate its equivalent form, which we will actually prove.

By an algebra we will always mean an associative commutative algebra with unit.

Recall that for any formal group law $F$ over a $\Q$-algebra $R$ there exists a unique \emph{logarithm} $g_F$, i.e., a power series $g_F(t)=t+ \ldots \in R[[t]]$ such that $g_F(F(x,y))=g_F(x)+g_F(y)$. In other words, $g_F$ is a strict isomorphism from the formal group law $F$ to the additive formal group law $\mathbb G_a$. For a complex Hirzebruch genus $\varphi\colon \varOmega^U_*\to R$ we will denote by $g_\varphi$ the logarithm of the formal group law $\varphi(F_U)$, the image of the universal formal group law in complex cobordism under $\varphi$. We have $g_\varphi=\varphi(g_U)$, where $g_U$ is the logarithm of the formal group law $F_U$ (over $\varOmega^U_*\otimes\Q$).

On the other hand, it is well known that the universal formal group law $F_O$ of unoriented geometric cobordism also has a logarithm $g_O(u)=u+\ldots \in \varOmega^O_*[[u]]$ such that $g_O(F_O(u,v))=g_O(u)+g_O(v)$ (for example, because a formal group law $F$ over a $\Z_2$-algebra is strictly isomorphic to the additive one if and only if $F(u,u)=0$, and this holds for $F_O$). So for any Hirzebruch genus $\varphi\colon \varOmega^O_*\to R$ there exists a logarithm $g_\varphi=\varphi(g_O)$ of the corresponding formal group law $F_\varphi=\varphi(F_O)$. In contrast to the rational case, this logarithm is not uniquely defined, but only up to postcomposition with the strict automorphisms of the additive formal group law. But any such automorphism (over a $\Z_2$-algebra $R$) has the form $a(x)=x+\sum a_k x^{2^k}$, in particular, $a'(x)=1$, and hence $g'_\varphi(u)$ is uniquely defined.

In fact, for any formal group law $F$, the derivative $g_F'(u)$ of any logarithm is always uniquely determined and equals $\frac{1}{\frac{\partial F(u,v)}{\partial v}|_{v=0}}=\omega_F$, the so-called \emph{invariant differential } of the formal group law $F$. Details on formal group laws in the context of bordism can be found, for example, in \cite{rav} or in the survey \cite{buch12}.

\begin{theorem}\label{th}
\begin{enumerate}[label=(\alph*)]
   \item  Let $R$ be a $\Q$-algebra and $\varphi \colon \varOmega^U_* \to R$ be a complex Hirzebruch genus with logarithm $g_\varphi(t) \in R[[t]]$. Then the following equality of power series in $R[[t]]$ holds: $$\varphi(SP(C_g,t)):=1 + \sum_{n=1}^{\infty} \varphi([SP^n(C_g)])t^n = \Bigl(\frac{1}{g_\varphi'(t)}\Bigr)^{g-1}=(\omega_\varphi(t))^{1-g}$$

   \item Let $R$ be a $\Z_2$-algebra and $\varphi \colon \varOmega^O_* \to R$ be an unoriented Hirzebruch genus with a logarithm $g_\varphi(t) \in R[[t]]$. Then the power series $g'(t)$ is even, so $g'_\varphi(\sqrt{t}) \in R[[t]]$ is well defined, and the following equality of power series in $R[[t]]$ holds: $$\varphi(SP(N_g,t)):=1 + \sum_{n=1}^{\infty} \varphi([SP^n(N_g)])t^n = \Bigl(\frac{1}{g'_\varphi(\sqrt t)}\Bigr)^{g-1}=\bigl(\omega_\varphi(\sqrt{t})\bigr)^{1-g}$$
   \end{enumerate}
\end{theorem}

\begin{proof}[Proof of the equivalence Theorem \ref{th1} $\iff$ Theorem \ref{th}.]
Consider the universal genera $\varphi_{univ}^U\colon\varOmega^U_*\hookrightarrow \varOmega^U_*\otimes\Q$ and $\varphi_{univ}^O\colon \varOmega^O_* \xrightarrow{\rm id} \varOmega^O_*$ respectively.

In the first case, by Mishchenko's formula \cite[Appendix 1]{nov2}, we have that the universal logarithm $g_U$ has the form $$g_{U}(u)=u + \sum_{n \ge 1} \frac{[\C P^n]}{n+1}u^{n+1},$$ and hence $\omega_U(u)=g'_{U}(u)=\C P(u)$. So Theorem \ref{th1} a) is a special case of Theorem \ref{th} a), and Theorem \ref{th} a) is obtained from Theorem \ref{th1} a) by applying the genus $\varphi$, since $g_\varphi=\varphi(g_{U})$.

In the unoriented case, the same argument works because we have $g_{O}'(u)=\R P(u^2)$ (as we mentioned above, this derivative does not depend on the ambiguity in the choice of logarithm $g_O$). In particular, it is even. Indeed, we have the following universal Buchstaber formula for the formal group law of geometric cobordism $F_{O}(x,y)=\frac{\sum [H_{ij}]x^i y^j}{\R P(x) \R P(y)}$, where $H_{ij}\subset \R P^i\times \R P^j$ is a hypersurface dual to the tensor product of tautological line bundles (see, for example, \cite[Theorem E.2.3]{bp}, whose proof also works in the unoriented case). By definition, we have $H_{0i}=\R P^{i-1}$ and one can easily check that $[H_{1i}]=[\R P^1][\R P^{i-1}]=0$ in $\varOmega^O_*$ (for example, by computing Stiefel--Whitney numbers). Then a direct calculation yields $g_{O}'(u)=\frac{1}{\frac{\partial F_O}{\partial v}|_{v=0}}=1+\sum [\R P^n]u^n$. Since $[\R P^{2n+1}]=0$, we finally obtain $\omega_O(u)=g'_{F_O}(u)=1+\sum [\R P^{2n}]u^{2n}=\R P(u^2)$.
\end{proof}

\begin{remark}
    Clearly for any power series $g(u)$ over a $\Z_2$-algebra $R$ its derivative $g'$ is always even.
\end{remark}

We also have the following equivalent form of Theorem \ref{th1} a).

\begin{corollary}\label{b}
If compact complex curves $C$ and $C'$ (not necessarily connected) are complex bordant, then $SP^n(C)$ and $SP^n(C')$ are also complex bordant. In particular, if a complex curve $C$ is complex null-bordant, then $SP^n(C)$ is also complex null-bordant.
\end{corollary}

\begin{proof}[Proof of the equivalence Theorem \ref{th1} a) $\iff$ Corollary \ref{b}]
    Let $C=\bigsqcup C_{g_i}$ and $C'=\bigsqcup C_{g'_i}$. Since complex bordism classes in $\varOmega^U_2$ are uniquely defined by the $c_1$ Chern number, we have $[C_g]=(1-g)[\C P^1]$. Since $[\C P^1]$ is a generator of the group $\varOmega_2^U\cong\Z$, $[C]=\sum_i(1-g_i)[\C P^1]$ is equal to $[C']=\sum_i(1-g'_i)[\C P^1]$ $\iff$ $\sum(g_i-1)=\sum(g'_i-1)$. But in this case, by Theorem \ref{th1} and formula \eqref{sum} we obtain \begin{multline*}
        SP(C,t)=\prod SP(C_{g_i},t)=\prod\frac{1}{\C P(t)^{g_i-1}}=\frac{1}{\C P(t)^{\sum (g_i-1)}}=\\
        =\frac{1}{\C P(t)^{\sum (g'_i-1)}}=\prod\frac{1}{\C P(t)^{g'_i-1}}=\prod SP(C_{g'_i},t)=SP(C',t)
        \end{multline*}

        On the other hand, if $SP$ is invariant under complex bordism, then from the equality $[C_g\sqcup \bigsqcup_{i=1}^{g-1}\C P^1]=0$ in $\varOmega^U_2$, we obtain $1=SP(C_g\sqcup \bigsqcup_{i=1}^{g-1}\C P^1,t)=SP(C_g,t)\C P(t)^{g-1}$, which is equivalent to Theorem \ref{th1} a).
\end{proof}

\begin{remark}
    An analogue of Corollary \ref{b} is completely false in the unoriented case. In fact, in $\varOmega^O_2=\Z_2$ we have $[N_{2k+1}]=0$ and $[N_{2k}]=[\R P^2]$, so the bordism class of $N_g$ depends only on the parity of $g$.
    
    On the other hand, all power series $\R P(t)^{k}$ are distinct, in other words, the bordism classes of all $SP^n(N_g)$ determine the integer $g$ itself. In particular, $\R P(t)^k=1$ only for $k=0$; i.e., $SP^n(N)$ is null-bordant for all $n$ only for the Klein bottle $N_1$.
\end{remark}

\begin{remark}
 It is not true that if an orientable surface $C$ is unoriented null-bordant, then $SP^n(C)$ is also unoriented null-bordant. For example, $[S^2]=0$ in $\varOmega^O_2$, but $[SP^{2n}(S^2)]=[\C P^{2n}]\ne 0$ in $\varOmega^O_{4n}$.
\end{remark}

\begin{remark}
    Complex bordism invariance is not true for symmetric products in general. For example, for zero-dimensional complex bordism, we have that $\pt_+ \sqcup\pt_-$ is null-bordant, but $SP^2(\pt_+\sqcup \pt_-)=\pt_+\times \pt_+ \sqcup \pt_+\times \pt_- \sqcup \pt_-\times \pt_- = \pt_+ \sqcup \pt_+ \sqcup \pt_-$ is not null-bordant.

It would be interesting to understand the geometric nature of the bordism invariance in Corollary \ref{b}.
\end{remark}

\begin{remark}
    In the paper \cite{egl} the authors proved that the complex bordism class of a Hilbert scheme of points of a smooth complex surface $S$ depends only on the bordism class of $S$ (by Fogarty's theorem \cite{fog} the Hilbert scheme of a smooth surface is smooth; this is not true in higher dimensions). More precisely, by algebro-geometric considerations they proved that the Chern numbers of ${\rm Hilb}^n(S)$ can be expressed as polynomials in the Chern numbers of $S$. Since for a complex curve the Hilbert scheme coincides with the symmetric product, we think that it is more or less clear that the arguments of \cite{egl} work (perhaps more easily) for complex curves and give Corollary \ref{b}, which is equivalent to Theorem \ref{th} a) formally.

    On the other hand, our arguments are completely different from those of \cite{egl}. They are more topological (and more classical in the context of symmetric products, we believe) and are based on the careful study of tangent bundles and cohomology rings of symmetric products. Of course, the Chern-number calculation is also hidden behind our bordism-theoretic considerations. But the bordism interpretation allows us to obtain Theorem \ref{th} first, which gives explicit formulas for bordism classes of symmetric products (not an abstract bordism-invariance result via implicit polynomial expressions in Chern numbers). Finally, this approach works in the unoriented case, where there is no bordism invariance at all.
\end{remark}

\section{Preliminaries}

For a genus $g$ curve $C_g$ and a point $x_0\in C_g$ consider the corresponding map $$\pi=\pi_{x_0}\colon SP^n(C_g)\to J(C_g), \hspace{4ex} [x_1, \ldots, x_n] \mapsto \sum x_i-nx_0,$$ where $J(C_g)$ is the Jacobian of the curve $C_g$ (considered as the variety of zero degree divisors on $C_g$). The Jacobian $J(C_g)$ is a complex torus, so topologically it is $T^{2g}$.

We have the following result of Mattuck and Macdonald.

\begin{fact}[\cite{mat1}, \cite{mcd}]
    If $n\ge 2g$, then $\pi\colon SP^n(C_g)\to T^{2g}$ is the projectivization of a holomorphic vector bundle $V^{n-g+1}$ (of complex dimension $n-g+1$).

    Moreover, there exist elements $x_i, x'_i\in H^1(T^{2g},\Z)$, $1\le i\le g$, such that $$H^*(T^{2g};\Z)=\Lambda_\Z(x_i,x'_i)$$ (exterior algebra over integers) and, for the total Chern class, $$c(V^{n-g+1})=\prod_{i=1}^g (1-x_ix'_i)$$
\end{fact}

If we denote by $\eta$ the tautological line bundle over $SP^n(C_g)=\mathbb P(V^{n-g+1})$ and $v=c_1(\bar\eta)$, then $c_{n-g+1}(\bar\eta\otimes V^{n-g+1})=v^{n-2g+1}\prod_{i=1}^g(v-x_ix'_i)$, and therefore $$H^*(SP^n(C_g);\Z)=\Lambda_\Z(x_i, x'_i)[v]/\bigl(v^{n-2g+1}\prod_{i=1}^g(v-x_ix'_i)=0\bigr), \quad n\ge 2g$$

Now consider the embedding $$j=j_{x_0}\colon SP^{n-1}(C_g)\hookrightarrow SP^n(C_g), \hspace{4ex} [x_1, \ldots, x_n]\mapsto [x_0, x_1, \ldots, x_n]$$ We will denote the restrictions $j^*\eta$ and $j^*(v)$ just by $\eta$ and $v$, respectively. If $n-1\ge 2g$, then $j$ is a projective subbundle (see \cite{mcd}), so this agrees with the previous notation. Similarly, the elements $x_i$, $x'_i$ are compatible with $j$, since $\pi\circ j=\pi$. So we have bundles $\eta$ and elements $v,x_i,x'_i$ for any $n$.

\begin{fact}[\cite{sch}]\label{sch}
    The embedding $j\colon SP^{n-1}(C_g)\hookrightarrow SP^n(C_g)$ is dual to the line bundle $\bar\eta$, i.e., the normal bundle $\nu(SP^{n-1}(C_g)\subset SP^n(C_g))$ is isomorphic to the restriction $\bar\eta|_{SP^{n-1}(C_g)}$, and for the fundamental class $[SP^{n-1}(C_g)] \in H_{2n-2}(SP^{n-1}(C_g);\Z)$ its $j_*$-image is Poincar\'e dual to the class $v\in H^2(SP^n(C_g);\Z)$.
\end{fact}

Now consider an antiholomorphic involution $\tau\colon C_g\to C_g$ without fixed points, such that $C_g/\tau = N_g$. It induces an involution $[x_1,\ldots, x_m] \mapsto [\tau (x_1),\ldots, \tau(x_m)]$ on $SP^{m}(C_g)$ for each $m$. Since $\tau$ acts freely on $C_g$, one can easily check that the fixed-point set $SP^m(C_g)^\tau$ can be identified with $SP^{m/2}(N_g)$ when $m$ is even and
is empty when $m$ is odd (see Remark \ref{taufix}).

 More precisely, we have the following result of Dupont and Lusztig.

\begin{fact}[\cite{dulu}]\label{dl}
    Let $n\ge g$. Then $\tau$ can be naturally extended to the projective bundle $\pi\colon SP^{2n}(C_g)\to T^{2g}$ and its fixed-point set can be identified with the real projectivization $SP^n(N_g)\to T^g$ of some real vector bundle $W^{2n-g+1}$.

    Moreover, there exist elements $y_i, y'_i\in H^1(T^{g},\Z_2)$, $1\le i\le \lfloor \frac{g}{2}\rfloor$, and, if $g$ is odd, an additional element $z\in H^1(T^g;\Z_2)$, such that $H^*(T^{g};\Z_2)=\Lambda_{\Z_2}(y_i, y'_i,z)$ and $$w(W^{2n-g+1})=\begin{cases}
        (1+y_1y'_1)\cdots (1+y_ky'_k),& g=2k\\
        (1+y_1y_1')\cdots (1+y_ky'_k)(1+z), & g=2k+1
    \end{cases}$$
\end{fact}

If we denote by $\xi$ the tautological line bundle over $SP^n(N_g)=\mathbb P(W^{2n-g+1})$ and $w=w_1(\xi)$, then $$w_{2n-g+1}(\xi\otimes W^{2n-g+1})=\begin{cases}
    w^{2(n-g)+1}(w^2+y_1y'_1)\cdots (w^2+y_ky'_k),& g=2k\\
        w^{2(n-g)+1}(w^2+y_1y'_1)\cdots (w^2+y_ky'_k)(w+z), & g=2k+1
\end{cases}$$ and $$H^*(SP^n(N_g);\Z_2)=\Lambda_{\Z_2}(y_i)[w]/\bigl(w_{2n-g+1}(\xi\otimes W^{2n-g+1})=0\bigr), \quad n\ge g$$

We also have similar embeddings $j\colon SP^{n-1}(N_g)\hookrightarrow SP^n(N_g)$. Similarly, we will denote $j^*\xi$, $j^*(w)$, and $j^*(y_i)$ simply by $\xi$, $w$, and $y_i$, respectively. Again, this notation agrees with the previous one, so we have bundles $\xi$ and elements $w, y_i$ for any $n$ (and $z$ for odd $g $).

From Theorems \ref{sch} and \ref{dl} we obtain the following

\begin{proposition}\label{dual}
    The embedding $j\colon SP^{n-1}(N_g)\hookrightarrow SP^n(N_g)$ is dual to the line bundle $\xi\oplus\xi$, i.e., the normal bundle $\nu(SP^{n-1}(N_g)\subset SP^n(N_g))$ is isomorphic to the restriction $(\xi\oplus\xi)|_{SP^{n-1}(N_g)}$, and for the fundamental class $[SP^{n-1}(N_g)] \in H_{2n-2}(SP^{n-1}(N_g);\Z_2)$ its $j_*$-image is Poincar\'e dual to the class $w^2\in H^2(SP^n(N_g);\Z_2)$.
\end{proposition}

We also have the following general results of Steenrod--Dold and Nakaoka.

\begin{fact}[\cite{dold}, \cite{nak}]
    For any finite CW complex $X$ and any abelian group $G$, the embedding $j\colon SP^{n-1}(X)\hookrightarrow SP^n(X)$ induces (splitting) monomorphisms on homology $$j_*\colon H_*(SP^{n-1}(X);G)\hookrightarrow H_*(SP^n(X);G)$$ and (splitting) epimorphisms on cohomology $$j^*\colon H^*(SP^{n}(X);G)\twoheadrightarrow H^*(SP^{n-1}(X);G)$$
\end{fact}

Hence, in our situation, we have (splitting) injections
\begin{gather*}
    j_*\colon H_*(SP^{n-1}(C_g);\Z)\hookrightarrow H_*(SP^n(C_g);\Z)\\
    j_*\colon H_*(SP^{n-1}(N_g);\Z_2)\hookrightarrow H_*(SP^n(N_g);\Z_2)
\end{gather*}
 and (splitting) surjections
  \begin{gather*}
    j^*\colon H^*(SP^{n}(C_g);\Z)\twoheadrightarrow H^*(SP^{n-1}(C_g);\Z)\\
    j^*\colon H^*(SP^{n}(N_g);\Z_2)\twoheadrightarrow H^*(SP^{n-1}(N_g);\Z_2)
\end{gather*}

In particular, we see that the cohomology rings $H^*(SP^n(C_g);\Z)$ are generated by the elements $x_i$, $x'_i$ and $v$ for any $n$ (since this is true for $n \ge 2g$). Similarly, the cohomology rings $H^*(SP^n(N_g);\Z_2)$ are generated by the elements $y_i$, $y_i'$, $z$ and $w$.

Since our symmetric products are oriented (with respect to the appropriate coefficients), we have the Poincar\'e-dual injections on cohomology
\begin{gather*}
    j_!\colon H^*(SP^{n-1}(C_g);\Z)\hookrightarrow H^{*+2}(SP^n(C_g);\Z)\\
    j_!\colon H^*(SP^{n-1}(N_g);\Z_2)\hookrightarrow H^{*+2}(SP^n(N_g);\Z_2)
\end{gather*}

We have the following easy consequence of the projection formula.

\begin{lemma}\label{l1}
    Let $\alpha\in H^*(SP^{n-1}(C_g);\Z)$. Consider any element $\alpha'\in H^*(SP^n(C_g);\Z)$ such that $j^*(\alpha')=\alpha$ (e.g., the same polynomial expression in the variables $x_i, x'_i, v$ as for the element $\alpha$, but considered as an element in $H^*(SP^n(C_g);\Z)$). Then $j_!(\alpha)=v\alpha'$.

    Similarly, for $\alpha\in H^*(SP^{n-1}(N_g);\Z_2)$ and the corresponding $\alpha'\in H^*(SP^n(N_g);\Z_2)$ we have $j_!(\alpha)=w^2\alpha'$.
\end{lemma}

\begin{proof}
    By definition $j_!=D^{-1}_nj_*D_{n-1}$, where $D_n$ and $D_{n-1}$ are the Poincar\'e duality operators in $SP^n$ and $SP^{n-1}$ respectively. So $j_!(\alpha)=D_n^{-1}j_*(\alpha\cap [SP^{n-1}])=D_n^{-1}j_*(j^*(\alpha')\cap [SP^{n-1}])$. Then, by the projection formula, we obtain $j_!(\alpha)=D_n^{-1}(\alpha'\cap j_*[SP^{n-1}])$. By Theorem \ref{sch} and Proposition \ref{dual} we have that $j_*[SP^{n-1}]$ is $D_n$-dual to $v$ or $w^2$ respectively. Denote this element uniformly by $u$. Then $j_!(\alpha)=D^{-1}_n(\alpha'\cap(u\cap [SP^n]))=D_n^{-1}((\alpha'u)\cap [SP^n])=\alpha'u$.
\end{proof}

\section{Proof of the main results}

We do not need the full cohomology ring structure of symmetric products, only the following relations in the top cohomology (because we are only interested in the calculation of characteristic numbers).

\begin{lemma}\label{l2}
    Denote by $s_i=x_ix'_i\in H^2(SP^n(C_g);\Z)$. Then, for any monomial $s_I$ with $|I|\le n$, we have the following equality in the top cohomology $H^{2n}(SP^n(C_g);\Z)$: $$v^{n-|I|}s_I=v^n$$

    Similarly, if we denote by $s_i=y_iy'_i\in H^2(SP^n(N_g);\Z_2)$, then for any monomial $s_I$ with $|I|\le n$, we have the following equality in the top cohomology $H^{2n}(SP^n(N_g);\Z_2)$: $$w^{2(n-|I|)}s_I=w^{2n}$$ And additionally, if $g=2k+1$, $$w^{2(n-|I|)-1}s_Iz=w^{2n}$$
\end{lemma}

\begin{proof}
    Let us prove the unoriented case for $g=2k$. The complex and $g=2k+1$ cases are similar.

    \emph{Case $n\ge g$.} We have the basic relation $$w^{2(n-g)+1}(w^2+s_1)\cdots(w^2+s_k)=0\quad \mbox{in }H^{2n-g+1}(SP^n(N_g);\Z_2)$$

    Multiplying by $\prod\limits_{j\notin I}(w^2+s_j)$ and using $s_i^2=0$ we obtain $w^{2n+1-4|I|}\prod\limits_{i\in I}(w^2+s_j)=0$, and hence in the top cohomology $w^{2n-2|I|}\prod\limits_{i\in I}(w^2+s_j)=0$. We can assume by induction on $|I|$ that $w^{2(n-|J|)}s_J=w^{2n}$ for $|J|<|I|$ (with the empty monomial as the base case). Then, expanding brackets, we obtain the required equation $w^{2(n-|I|)}s_I=w^{2n}$.

    \emph{Case $n<g$.} By descending induction on $n$, assume that we know the required equation in the top cohomology of $SP^n(N_g)$. Then, by Lemma \ref{l1}, we have the injection $j_!\colon H^{2n-2}(SP^{n-1}(N_g);\Z_2)\hookrightarrow H^{2n}(SP^n(N_g);\Z_2)$, $\alpha\mapsto w^2\alpha'$. For $\alpha=w^{2n-2-2|I|}s_I-w^{2n-2}$ we have that $\alpha'=w^{2n-2-2|I|}s_I-w^{2n-2}$, and hence $j_!(\alpha)=w^{2n-2|I|}s_I-w^{2n}=0$. Since $j_!$ is injective, we obtain $\alpha=0$, as claimed.
\end{proof}

\begin{remark}
   In the complex case, these relations are due to Macdonald \cite{mcd}.\end{remark}
    
   \begin{remark}
        A complete description of the cohomology ring $H^*(SP^n(C_g);\Z)$ can be found in \cite{gug4}.
\end{remark}

\begin{proof}[Proof of Theorem \ref{th}]

\emph{Complex case}

For $n\ge 2g$ we have $SP^n(C_g)=\mathbb P(V^{n-g+1})$, so the total Chern class has the form
\begin{equation}\label{chern}
    c(SP^n(C_g))=c(\bar\eta\otimes V^{n-g+1})=(1+v)^{n-2g+1}\prod_{i=1}^g(1+v-s_i)
\end{equation}

Since by Theorem \ref{sch} we have that $$c(SP^{n-1}(C_g))=j^*\Bigl(\frac{c(SP^n(C_g))}{c(\bar\eta)}\Bigr)$$ we obtain that formula \eqref{chern} in fact holds for any $n$.

Now let $\varphi\colon\varOmega^U_*\to R$ be a Hirzebruch genus with exponential $f_\varphi(t)$ and logarithm $g_\varphi(u)$. Then by the Hirzebruch formula  $$\varphi([SP^n(C_g)])=\langle\prod_{i=1}^n \frac{t_i}{f_\varphi(t_i)}, [SP^n(C_g)]_\Z\rangle,$$
where $t_i$ are the Chern roots of the tangent bundle $TSP^n(C_g)$, and $[SP^n(C_g)]_\Z$ is the integer fundamental class.

The Chern roots can be easily seen from \eqref{chern}, so we obtain  \begin{equation}\label{hirz}
    \varphi([SP^n(C_g)])=\langle\Bigl(\frac{v}{f_\varphi(v)}\Bigr)^{n-2g+1}\prod_{i=1}^g \frac{v-s_i}{f_\varphi(v-s_i)}, [SP^n(C_g)]_\Z\rangle
\end{equation}

Since $s_i^2=0$ we have
    \begin{equation}\label{fr}\frac{v-s_i}{f_\varphi(v-s_i)}=\frac{v}{f_\varphi(v)}+s_i \cdot \frac{\partial}{\partial s_i}\Bigl(\frac{v-s_i}{f_\varphi(v-s_i)}\Bigr)\bigg|_{s_i=0}=\frac{v}{f_\varphi(v)}+s_i \frac{-f_\varphi(v)+f'_\varphi(v)v}{f_\varphi(v)^2}\end{equation}

    Now we need to consider the $2n$-th graded component of the right-hand side of \eqref{hirz} and evaluate it on the fundamental class.

    Since by Lemma \ref{l2} in top degree we have $s_I v^k =v^n$, we can replace all $s_i$ in \eqref{fr} by $v$ and obtain $\bigl(\frac{v}{f_\varphi(v)}\bigr)^2f'_\varphi(v)$.

Therefore, finally
\begin{multline*}
    \varphi([SP^n(C_g)])=\langle\bigl(f'_\varphi(v)\bigr)^g\Bigl(\frac{v}{f_\varphi(v)}\Bigr)^{n+1}, [SP^n(C_g)]_\Z\rangle={\rm coeff}_n \bigl(f'_\varphi(v)\bigr)^g\Bigl(\frac{v}{f_\varphi(v)}\Bigr)^{n+1}=\\
    =\frac{1}{2\pi i} \oint \frac{(f'_\varphi(t))^g}{f_\varphi(t)^{n+1}} dt=\bigg| {\substack{u=f_\varphi(t),\, t=g_\varphi(u)\\ dt=g'_\varphi(u)du=\frac{du}{f'_\varphi(t)}}} \bigg|= \frac{1}{2\pi i} \oint \frac{1}{g'_\varphi(u)^{g-1}u^{n+1}} du= {\rm coeff}_n \frac{1}{(g'_\varphi(u))^{g-1}}
\end{multline*}

   which is equivalent to the statement of Theorem \ref{th}.

 \vspace{1ex}

   \emph{Unoriented case}

The unoriented case can be proved in two ways.

\emph {First proof.} This is the same argument as in the complex case. Consider, for example, the case $g=2k+1$. The even genus case is similar.

  Again, we have the total Stiefel--Whitney class of the following form
\begin{equation}\label{whit}
    w(SP^n(N_g))=(1+w)^{2n-2g+1}(1+w+z)\prod_{i=1}^k(1+w^2+s_i)
\end{equation}

Now let $\varphi\colon\varOmega^O_*\to R$ be an unoriented Hirzebruch genus with exponential $f_\varphi(t)$ and logarithm $g_\varphi(u)$. Then  $$\varphi([SP^n(N_g)])=\langle\prod_{i=1}^{2n} \frac{t_i}{f_\varphi(t_i)}, [SP^n(N_g)]_{\Z_2}\rangle,$$
where $t_i$ are the Stiefel--Whitney roots of $TSP^n(N_g)$ (see the first proof of Theorem E.3.3 in \cite{bp}, which also works in the unoriented case).

From \eqref{whit} we obtain \begin{equation}\label{hirz2}
    \varphi([SP^n(N_g)])=\langle\Bigl(\frac{w}{f_\varphi(w)}\Bigr)^{2n-2g+1}\frac{w+z}{f_\varphi(w+z)}\prod_{i=1}^k \frac{(w+t_i)^2}{f_\varphi(w+t_i)^2}, [SP^n(N_g)]_{\Z_2}\rangle
\end{equation}
where $t_i^2=s_i$.

Since $t_i^4=s_i^2=0$ we have
    \begin{equation}\label{fr2}\Bigl(\frac{w+t_i}{f_\varphi(w+t_i)}\Bigr)^2=\Bigl(\frac{w}{f_\varphi(w)}+t_i \frac{f_\varphi(w)+f'_\varphi(w)w}{f_\varphi(w)^2}\Bigr)^2=\frac{w^2}{f_\varphi(w)^2}+s_i \frac{f_\varphi(w)^2+f'_\varphi(w)^2w^2}{f_\varphi(w)^4}\end{equation}

    Similarly, \begin{equation}\label{fr3}\frac{w+z}{f_\varphi(w+z)}=\frac{w}{f_\varphi(w)}+z \frac{f_\varphi(w)+f'_\varphi(w)w}{f_\varphi(w)^2}\end{equation}

    As we need to consider the $2n$-th graded component of the right-hand side of \eqref{hirz2} and evaluate it on the fundamental class, and since by Lemma \ref{l2} in top-degree cohomology we have $s_I w^k =s_Iw^{k-1}z =w^{2n}$, we can replace all $s_i$ in \eqref{fr2} by $w^2$ and $z$ in \eqref{fr3} by $w$ and obtain $\bigl(\frac{w^2}{f_\varphi(w)^2}\bigr)^2f'_\varphi(w)^2$ and $(\frac{w}{f_\varphi(w)})^2f'_\varphi(w)$ respectively.

Therefore, finally
\begin{multline*}
    \varphi([SP^n(N_g)])=\langle\bigl(f'_\varphi(w)\bigr)^g\Bigl(\frac{w}{f_\varphi(w)}\Bigr)^{2n+1}, [SP^n(N_g)]_\Z\rangle={\rm coeff}_{2n} \frac{w^{2n+1}f'_\varphi(w)^g}{f_\varphi(w)^{2n+1}}=\\
    =\mathrm{res}_0\Bigl(\frac{f'_\varphi(t)^g}{f_\varphi(t)^{2n+1}}\Bigr)= \bigg| {\substack{u=f_\varphi(t),\, t=g_\varphi(u)\\ dt=g'_\varphi(u)du=\frac{du}{f'_\varphi(t)}}} \bigg|= \mathrm{res}_0\Bigl(\frac{1}{g'_\varphi(u)^{g-1}u^{2n+1}}\Bigr)= {\rm coeff}_{2n} \frac{1}{(g'_\varphi(u))^{g-1}}
\end{multline*}

which is equivalent to the required identity $\varphi(SP(N_g,u))=\bigl(\frac{1}{g'_\varphi(\sqrt{u})}\bigr)^{g-1}$.

\smallskip
\emph{Second proof.}

In fact, the formula for $SP^n(N_g)$ follows from the formula for $SP^n(C_g)$.
Choose an antiholomorphic involution $\tau$ on $C_g$ with $C_g/\tau=N_g$.
As was mentioned above, $\tau$ induces an involution on $SP^{m}(C_g)$ for each $m$, and the fixed-point set is $SP^{m/2}(N_g)$ when $m$ is even, and
is empty when $m$ is odd.

We have the following general result of Conner and Floyd.

\begin{fact}[{\cite[Theorem 24.4]{cf}}]
    If an antiholomorphic involution on an almost complex manifold $M$ has fixed-point set $F$, then $r_*[M]=[F]^2$, where $r_*\colon\varOmega_*^U\to \varOmega^O_*$ is the forgetful map. 
\end{fact}

Applying this
to $\C P^n$ with complex conjugation gives $r_*([\C P^n])=[\R P^n]^2$, hence
\[
r_*(\C P(t))=1+\sum_{n\ge 1}[\R P^n]^2t^n
=1+\sum_{n\ge 1}[\R P^{2n}]^2t^{2n}=\R P(t)^2,
\]
where we used $[\R P^{2n+1}]=0$ and the fact that, modulo $2$,
$\R P(t)^2=1+\sum_{n\ge 1}[\R P^{2n}]^2t^{2n}$.

On the other hand, applying the Conner--Floyd theorem to the extended involution on $SP^{2n}(C_g)$ gives
$r_*([SP^{2n}(C_g)])=[SP^n(N_g)]^2$ and $r_*([SP^{2n+1}(C_g)])=0$. Therefore
\[
r_*(SP(C_g,t))
=\sum_{n\ge 0} r_*([SP^{2n}(C_g)])t^{2n}
=\sum_{n\ge 0}[SP^n(N_g)]^2t^{2n}
=SP(N_g,t)^2.
\]

On the other hand, by Theorem \ref{th1} a) and the multiplicativity of $r_*$,
\[
r_*(SP(C_g,t))=r_*\left(\frac{1}{\C P(t)^{g-1}}\right)
=\frac{1}{r_*(\C P(t))^{g-1}}
=\frac{1}{\R P(t)^{2(g-1)}}.
\]
Combining the two computations,
\[
SP(N_g,t)^2=\frac{1}{\R P(t)^{2(g-1)}}.
\]
Finally, since the squaring map
$\varOmega^O_*\xrightarrow{x\mapsto x^2}\varOmega^O_*$ is injective, we obtain
\[
SP(N_g,t)=\frac{1}{\R P(t)^{g-1}}.
\]
\end{proof}

\section{Blagojevi\'{c}--Gruji\'{c}--Zivaljevi\'{c} conjecture for non-orientable surfaces}

It is well known that for any connected CW complex $X$ the inclusion $X \hookrightarrow SP^n(X)$ induces an abelianization map $\pi_1(X) \to \pi_1(SP^n(X))=\pi_1(X)^{ab}=H_1(X;\Z)$ on fundamental groups. In particular, we have an isomorphism on the first homology $H_1(X;\Z) \cong H_1(SP^n(X);\Z)$ (for a generalization to quotients of $X^{\times n}$ by any subgroup of $\Sigma_n$, see \cite{kt}).

Therefore, since the first homology distinguishes closed surfaces (orientable or not), we conclude that the first homology of any symmetric product of a closed surface also allows us to reconstruct the surface itself.

In particular, we obtain the following proposition.
\begin{proposition}
    If for two closed surfaces $S$ and $S'$ we have a homotopy equivalence $SP^n(S)\simeq SP^n(S')$ for some $n$, then in fact $S$ and $S'$ are homeomorphic (and hence, $SP^n(S)$ and $SP^n(S')$ are also homeomorphic).
\end{proposition}

On the other hand, let us consider surfaces of finite type, i.e., homeomorphic to a closed surface with finitely many punctures. Denote by $C_{g,s}$ a complex curve of genus $g$ with $s\ge 1$ punctures, and by $N_{g,s}$ a surface $N_g=(\R P^2)^{\#(g+1)}$ with $s\ge 1$ punctures. Then $C_{g,s}\simeq \bigvee_{i=1}^{2g+s-1} S^1$ and $N_{g,s}\simeq \bigvee_{i=1}^{g+s} S^1$. Since the symmetric product is homotopy invariant (and the number of circles is again determined by the fundamental group of the symmetric products), we obtain the following

\begin{proposition}
    We have $SP^n(N_{g,s}) \simeq SP^n(N_{g',s'})\iff N_{g,s}\simeq N_{g',s'}\iff g+s=g'+s'$.
    
    And similarly, $SP^n(C_{g,s}) \simeq SP^n(C_{g',s'})\iff C_{g,s}\simeq C_{g',s'}\iff 2g+s=2g'+s'$.
\end{proposition}

In fact, $SP^n(\bigvee_{i=1}^m S^1)\simeq {\rm sk}^n T^m$ (see \cite{ong}). So, in the noncompact finite-type case, the homotopy type of symmetric products carries exactly the same information as the homotopy type of the surface itself. It turns out that the topological type of any symmetric product allows us to reconstruct the topological type of the surface itself.

\begin{theorem}[non-orientable extension of Blagojevi\'{c}--Gruji\'{c}--Zivaljevi\'{c} conjecture \cite{bgz1, bgz2}]\label{BGZ}

If for two surfaces $S$ and $S'$ of finite type the symmetric products $SP^n(S)$ and $SP^n(S')$ are homeomorphic, then $S$ and $S'$ are also homeomorphic.

\end{theorem}

The case of oriented surfaces was proved by the second author in \cite{gug4}. Here we prove the remaining cases, following essentially the same approach.

\begin{remark}
    The following question remains: does the Blagojevi\'{c}--Gruji\'{c}--Zivaljevi\'{c} conjecture hold for infinite-type surfaces?
\end{remark}

\begin{proof}

Recall that for $n\ge g$ we have a projective bundle $\pi\colon SP^n(N_g)\to T^g$. From Theorem \ref{dl} we obtain that $$w(SP^n(N_g))=\begin{cases}
    (1+w)^{2(n-g)+1}\prod_{i=1}^k(1+w^2+y_iy'_i), & g=2k\\
    (1+w)^{2(n-g)+1}(1+w+z)\prod_{i=1}^k(1+w^2+y_iy'_i), & g=2k+1
\end{cases}$$
(this is true for the projective bundle, so for $n\ge g$; for smaller $n$ it follows, since the inclusion $j\colon SP^{n-1}(N_g)\to SP^{n}(N_g)$ is dual to $w^2$).

We only need the first two Stiefel--Whitney classes: $$w_1(SP^n(N_g))=\begin{cases}
    w, & g=2k\\
    z, & g=2k+1
\end{cases}$$

$$w_2(SP^n(N_g))=\begin{cases}
    (n-k)w^2+\sum_{i=1}^ky_iy'_i, & g=2k\\
    (n-k)w^2+wz+\sum_{i=1}^ky_iy'_i, & g=2k+1
\end{cases}$$

We have that the inclusion $N_g \hookrightarrow SP^n(N_g)$ induces an isomorphism $H^1(SP^n(N_g);\Z_2)\cong H^1(N_g;\Z_2)=(\Z_2)^{g+1}$, and the elements $y_i, y_i', z, w$ form a basis in this space.

Now, the inclusion $N_{g,s}\hookrightarrow N_g$ induces a monomorphism $H^1(N_g;\Z_2)\hookrightarrow H^1(N_{g,s};\Z_2)$, so we may assume that $y_i, y_i', w,z$ are linearly independent elements in $H^1(SP^n(N_{g,s});\Z_2)=V$. Recall that a symmetric product of any surface has a canonical smooth structure, natural with respect to inclusions of open subsets (see Proposition \ref{exist} and Remark \ref{BGZrem} in the Appendix). Hence the tangent bundle of $SP^n(N_{g,s})$ is just the restriction of the tangent bundle of $SP^n(N_g)$. In particular, we see that $w_1(SP^n(N_{g,s}))\ne 0$; in other words, $SP^n(N_{g,s})$ is nonorientable (this can also be seen in another way: for a non-orientable surface $S$, considering an embedded M\"obius band $M\hookrightarrow S$ and $n$ disjoint disks in its complement, we obtain an embedding of $M\times D^{2n}$ into the symmetric product $SP^{n+1}(S)$).

Since for an orientable surface its symmetric products are complex manifolds (see the Appendix), in particular, they are orientable. Since orientability is a topological invariant, it is sufficient to check that the topological type of a symmetric product of a nonorientable surface allows us to reconstruct the surface.

Since $SP^n(N_{g,s})\simeq SP^n(\bigvee _{i=1}^{g+s} S^1)\simeq {\rm sk}^n T^{g+s}$, we have that its cohomology is a truncated exterior algebra $\Lambda_{\Z_2}^{\le n}(V)$. We may assume that $n\ge 2$, since if $n=1$ the conjecture is trivial. Since in an exterior algebra the square of any one-dimensional element is zero, we obtain that \begin{equation*}
w_2(SP^n(N_{g,s}))=\begin{cases}
    \sum_{i=1}^k y_i y_i' & g=2k\\
    wz+\sum_{i=1}^{k-1} y_i y_i' & g=2k-1
\end{cases}
\end{equation*}

Let $SP^n(N_{g,s})\cong SP^n(N_{g',s'})$. Thus, we have an induced isomorphism $\Lambda_{\Z_2}^{\le n}(V)\cong \Lambda_{\Z_2}^{\le n}(V')$.

Now, notice that the Stiefel--Whitney classes are topological invariants; this was originally proved in \cite{thom}.

But one can easily see that if an isomorphism $\Lambda_{\Z_2}^{\le n}(V)\cong \Lambda_{\Z_2}^{\le n}(V')$ maps $\sum a_i b_i$ to $\sum a'_jb'_j$ (for linearly independent elements), then it must send a subspace $\langle a_i, b_i \rangle$ to a subspace $\langle a_j', b_j' \rangle$. In particular, they have the same dimension, and hence the number of summands in $\sum a_i b_i$ and $\sum a'_j b'_j$ is the same.

From this, we obtain that the only cases that remain to be distinguished are $N_{2k,s}$ and $N_{2k-1, s+1}$. In this case, the isomorphism sends $w_2=\sum_{i=1}^{k} y_i y_i'$ to $\tilde w_2=\tilde w z+ \sum_{i=1}^{k-1} \tilde y_i \tilde y'_i$, and hence maps $\langle y_i, y'_i\rangle$ isomorphically onto $\langle\tilde y_i, \tilde y_i', \tilde w, z\rangle$. But it also must send $w_1=w\notin \langle y_i, y_i'\rangle $ to $\tilde w_1=z\in \langle\tilde y_i, \tilde y_i', \tilde w, z\rangle$, which is a contradiction.
\end{proof}

\section{Appendix: smooth structures on symmetric products}

In this appendix, by a manifold we mean a second-countable Hausdorff topological manifold (without boundary). In particular, by a surface we mean a $2$-dimensional manifold.

Let $M$ be a manifold and $G\subset \Sigma_n$ be a subgroup. Let us start with the following question: when is $M^{\times n}/G$ a manifold? It is easy to see the following.

\begin{proposition}\label{nm}
    If $\dim M\ge 3$ and $G$ is nontrivial, then $M^{\times n}/G$ is not an integral homology manifold.
\end{proposition}

\begin{proof}
    Since the question is local, it is sufficient to check that $(\R^m)^n/G$ is not a homology manifold for $m\ge 3$.

   If $G$ is nontrivial, it contains an element $g$ of prime order $p$. In the disjoint cycle decomposition of $g \in G \subset \Sigma_n$ all cycles have order $p$. Choose any such $p$-cycle $(i_1, \ldots, i_p)$ permuted by $g$ and denote the corresponding coordinate subspace by $(\R^m)^I \subset (\R^m)^n$.

   Denote by $C_p$ the cyclic subgroup $\langle g \rangle$. The fixed-point set of this subgroup on the subspace $(\R^m)^I$ is the diagonal $D=\{(x, \ldots, x) \mid x \in \R^m\}$. The orthogonal complement $D^\perp \subset (\R^m)^I$ is also $g$-invariant and has dimension $(p-1)m\ge 3$.

   Now, pick a point $x \in (\R^m)^n$ with the stabilizer ${\rm Stab}_Q(x)=C_p$ (with $I$-coordinates in $D$ and with generic other coordinates). By the Slice Theorem (see, for example, \cite[Theorem B.24]{ggk}) the $G$-manifold $(\R^m)^n$ looks $G$-equivariantly locally near the orbit $G\cdot x$ like $G\times_{C_p}((\R^m)^n)$ near $G\cdot 0$. Hence, the quotient $(\R^m)^n/G$ looks locally near the image of $x$ like $(\R^m)^n/C_p$ near zero. Denote by $F$ the fixed-point subspace of $C_p$ in $(\R^m)^n$. It is easy to see that $C_p$ acts freely on the orthogonal complement $F^\perp$, and hence on the unit sphere $S(F^\perp)$. Then we obtain that the quotient $(\R^m)^n/G$ looks locally around zero like $(F\times F^\perp)/C_p=F \times (F^\perp/C_p)$.

    Now it is easy to see that $(\R^m)^n/G$ is not a homology manifold, since the link $F\times (F^\perp/C_p)-0$ is not an $(mn-1)$-dimensional homology sphere. Indeed, we have that $F^\perp/C_p$ is a cone over the lens space $S(F^\perp)/C_p$. Similarly, $F$ is a cone over the sphere $S(F)$. So $$F\times (F^\perp/C_p)-0=(F-0)\times (F^\perp/C_p)\cup_{(F-0)\times (F^\perp/C_p-0)} F\times (F^\perp/C_p-0)$$ is homotopy equivalent to the join $S(F)\ast (S(F^\perp)/C_p)\simeq \Sigma^{\dim F} S(F^\perp)/C_p$. Since $F^\perp \supset D^\perp$, $\dim F^\perp \ge 3$, and hence $S(F^\perp)$ is simply connected. Since $C_p$ acts freely on $S(F^\perp)$, we have $H_1(S(F^\perp)/C_p;\Z)\cong\pi_1(S(F^\perp)/C_p)\cong \Z/p$. Therefore, $H_{\dim F+1}(\Sigma^{\dim F} S(F^\perp)/C_p;\Z)\cong \Z/p$. But $\dim F+1<nm-1$.
\end{proof}

\begin{remark}
    We see that $M^{\times n}/G$ is not even a homology manifold with boundary.
\end{remark}

\begin{remark}
    The cases of the whole permutation group and of cyclic subgroups were proved in \cite{wa}.
\end{remark}

The case of $1$-dimensional manifolds is somewhat subtle. We have the following general theorem of Lange.

\begin{fact}[\cite{la}]
    For a finite subgroup $G\subset O(n)$ the quotient $\R^n/G$ is a (homology) manifold (with boundary) $\iff$ $G=G_{rr}\times P_1\times \ldots \times P_k$ acting on the orthogonal decomposition $\R^n=V_{rr}\oplus V_1\oplus\ldots\oplus V_k$, where $G_{rr} \subset O(V_{rr})$ is generated by reflections and rotations, and $P_i\subset SO(V_i)$ are conjugate to the Poincar\'e subgroup $P\subset SO(4)\subset SO(V_i)$ (and additionally there are no Poincar\'e factors if $n\le 4$, or $n=5$ and $G_{rr}$ contains a reflection). Moreover, $\R^n/G$ has no boundary $\iff$ $G$ does not contain reflections.
\end{fact}

In the case of a permutation subgroup $G\subset \Sigma_n$, it is easy to see that no Poincar\'e factor can occur, and we obtain the following.

\begin{corollary}
    For a $1$-dimensional manifold $M$ and a subgroup $G\subset \Sigma_n$, the quotient $M^{\times n}/G$ is a manifold with boundary $\iff$ $G$ is generated by transpositions, $3$-cycles and products of two disjoint transpositions; and the boundary is empty $\iff$ $G$ contains no transpositions.
\end{corollary}
\smallskip

Now let us consider the case of a surface $M$. Recall the notation $SP^n(M)=M^{\times n}/\Sigma_n$. Since $SP^n(\R^2)$ is homeomorphic to $\R^{2n}$, it is easy to see that the symmetric product $SP^n(M)$ is a manifold.

Let us assume that $G\subset \Sigma_n$ is generated by transpositions. This class of subgroups consists exactly of the Young subgroups $\Sigma_{I_1}\times\cdots\times \Sigma_{I_k}\subset \Sigma_n$ for partitions ${\bf n}:=\{1, \ldots, n\}=I_1\sqcup \ldots \sqcup I_k$. Thus, in this case we obviously have that $M^{\times n}/G \cong SP^{|I_1|}(M)\times\cdots\times SP^{|I_k|}(M)$ is a manifold. In fact, it is a smooth manifold, as we will see below.

    On the other hand, we have the following.

    \begin{proposition}
        If $G$ is not a Young subgroup, then $M^{\times n}/G$ is not an integral homology manifold.
    \end{proposition}

\begin{remark}
    This was essentially proved in \cite{ke01}. For the sake of completeness (and because the argument is short), we give a proof of this fact here.
\end{remark}

    \begin{proof}
        Since the question is local, it is sufficient to consider the case $M=\C$. The argument is similar to that in Proposition \ref{nm}.

        Consider the maximal Young subgroup $H \subset G$, i.e., generated by all transpositions in $G$. Then $H$ is $\Sigma_{I_1}\times\cdots\times \Sigma_{I_k}\subset G \subset \Sigma_n$, and it is normal in $G$. The blocks $I_j\subset {\bf n}$ are exactly the $H$-orbits. The group $G$ permutes these orbits, but preserves their cardinality. So the quotient group $Q=G/H$ acts on the set $\{I_1, \ldots, I_k\}$ and also preserves the cardinalities of the blocks (in particular, we have an embedding $Q \hookrightarrow \Sigma_k$). Notice that $Q$ does not contain (nontrivial) transpositions of blocks of cardinality $1$, since they correspond to transpositions in $G$, and hence are trivial in $Q$.

Since $G$ is not a Young subgroup, $Q$ is nontrivial, so it contains a nontrivial element $q\in Q$ of prime order $p$. Since any such element can be decomposed into disjoint cycles of order $p$ (in $\Sigma_k$), we can consider blocks $J_1,\ldots, J_p$ cyclically permuted by $q$ (under one of its prime cycles). Denote the common cardinality of $J_i$ by $m$. Notice that if $p=2$, then $m\ge 2$, since $Q$ does not contain transpositions of singleton blocks.

    We have $\C^n/H\cong SP^{|I_1|}(\C) \times \cdots \times SP^{|I_k|}(\C) \cong \C^n$. To distinguish the quotient from the source, let us denote the quotient by $\C^{\bf n}$. So the quotient group $Q$ acts on the quotient $\C^{\bf n}$ by permuting the $I$-blocks of coordinates and preserving the cardinalities of the blocks.

    Then the blocks $J_1, \ldots, J_p$ correspond to the coordinate subspace $(\C^m)^p \subset \C^{\bf n}$. The element $q$ preserves $(\C^m)^p$ and acts on it by cyclic permutation of the $J$-blocks. Hence, the fixed subspace of this restricted action is the diagonal subspace $D=\{(v,\ldots\, v) \mid v\in \C^m\} \subset (\C^m)^p$. Since the permutation action is obviously orthogonal, we have a $q$-invariant orthogonal complement $D^\perp \subset (\C^m)^p$ of complex dimension $(p-1)m\ge 2$.

Denote by $F$ the fixed-point subspace of $C_p$ in $\C^{\bf n}$. Since $F^\perp \supset D^\perp$, $\dim_\C F^\perp \ge 2$, and hence $S(F^\perp)$ is simply connected. Now, similarly to the proof of Proposition \ref{nm}, for a generic point $x \in \C^{\bf n}$ with the stabilizer ${\rm Stab}_Q(x)=C_p$, we obtain that the link of the image of $x$ in $\C^{\bf n}/Q=\C^n/G$ is homotopy equivalent to the join $S(F)\ast (S(F^\perp)/C_p)\simeq \Sigma^{\dim_\C F} S(F^\perp)/C_p$, which is not a homology sphere of dimension $2n-1$.
    \end{proof}

Now, for a surface $M$, we want to discuss smooth structures on $SP^n(M)$ (up to diffeomorphism). By Rad\'o's classical theorem \cite{rado}, every surface is triangulable (and hence has a smooth structure). In \cite{eps} it was proved that any homeomorphism between smooth surfaces is isotopic to a diffeomorphism (see also \cite{ha26}). Applying this result to the identity map on a surface equipped with two smooth structures, we obtain an isotopy (through homeomorphisms) between the identity and some diffeomorphism between the two smooth structures. Hence a smooth structure on a surface is canonical in the sense that it is unique up to diffeomorphism isotopic (through homeomorphisms) to identity. In particular, we obtain a canonical smooth structure on the Cartesian product $M^{\times n}$.

Moreover, since any orientable surface has a complex structure, there exists an atlas on it with only holomorphic transition functions. On the other hand, since any non-orientable surface can be obtained from a complex curve as the quotient by a free antiholomorphic involution, for any surface $M$ there exists a smooth atlas such that each transition function is either holomorphic or antiholomorphic (a so-called Klein, or dianalytic, structure on $M$; for details, see \cite{ag}).

Fix such a structure on $M$. Before we give a careful proof of the construction of a natural dianalytic structure on $SP^n(M)$ from that on $M$, we need the following well-known observation.

\begin{lemma}[Holomorphic Symmetric Function Theorem]
    Consider a branched covering map $\sigma\colon \C^n\to \C^n$, $(z_1, \ldots, z_n)\mapsto (\sigma_1(z), \ldots, \sigma_n(z))$, where $\sigma_i$ are the elementary symmetric polynomials.

    Let $U$ be a $\Sigma_n$-invariant open subset in $\C^n$.  Then $\sigma(U)$ is an open subset in $\C^n$, and $\sigma$ induces a bijection $\sigma^*\colon \mathcal O(\sigma(U))\to \mathcal O(U)^{\Sigma_n}$, $f\mapsto f\circ \sigma$, between holomorphic functions on $\sigma(U)$ and symmetric (i.e., $\Sigma_n$-invariant) holomorphic functions on $U$.
\end{lemma}

\begin{proof}
    The map $\sigma$ is a topological quotient map. In particular, $\sigma(U)$ is open, and for any continuous symmetric $f\colon U \to \C$ there exists a unique continuous $\tilde f\colon \sigma(U) \to \C$ such that $f=\tilde f \circ \sigma$. We need to check that $\tilde f$ is holomorphic if $f$ is.

Consider the polynomial $\prod_{i=1}^n (t-z_i)$. Its discriminant $\Delta$ is a symmetric polynomial in $(z_1, \ldots, z_n)$, hence can be written as a polynomial $P$ in $(\sigma_1, \ldots, \sigma_n)$. Denote by $D$ the divisor $P=0$. Then $\sigma$ is a regular holomorphic covering over the complement $\sigma(U)- D$. Therefore, $\tilde f$ is holomorphic on $\sigma(U)- D$.

But by the standard Riemann removable singularities theorem (see, for example, \cite[Chapter I, Section C, Theorem 3]{gr}), any locally bounded function $\tilde f$ holomorphic on $U-D$ is holomorphic on $U$ (for $n>1$).
\end{proof}

\begin{proposition}\label{exist}
\vphantom{g}
        There exists a functorial way to construct a dianalytic structure on $SP^n(M)$ from a Klein structure on $M$ (i.e., such that any dianalytic map of surfaces $M \to N$ induces a dianalytic map $SP^n(M) \to SP^n(N)$).

        If $M$ is complex, then we obtain a functorial complex structure on $SP^n(M)$.

        In particular, we obtain a smooth structure on $SP^n(M)$, natural with respect to dianalytic maps.
\end{proposition}

\begin{proof}
    1) \emph{Construction of an atlas}

Consider an arbitrary point $x=[x_1,\ldots, x_n]\in SP^n(M)$, $x_i\in M$. Assume that the distinct points among them are exactly $p_1,\ldots, p_k\in M$ with multiplicities $n_1,\ldots, n_k\ge 1$ respectively.

Choose pairwise disjoint open neighborhoods $U_j\ni p_j$ with diffeomorphisms (from the given Klein structure on $M$) $\varphi_i\colon U_i \to \mathbb D\subset \C$, where $\mathbb D$ is an open unit disk.

Consider an open subset $U = U_1^{\times n_1} \times U_2^{\times n_2} \times \cdots \times U_k^{\times n_k} \subset M^{\times n}$. Then $V = p(U)$ is an open neighborhood of $x$ in $SP^n(M)$. Notice that $U$ is invariant under the action of the subgroup $G=\Sigma_{n_1}\times\cdots\times \Sigma_{n_k}\subset \Sigma_n$. Since the $U_j$ are pairwise disjoint, one can easily check that the natural projection $p\colon M^{\times n} \to SP^n(M)$ induces a homeomorphism $\hat p\colon U/G\cong p(U)=V$.

The diffeomorphisms $\varphi_i\colon U_i \to \mathbb D$ induce a diffeomorphism
$\varphi=\varphi_1\times\ldots\times \varphi_k \colon U \xrightarrow{\cong} \mathbb D^{n_1} \times \cdots \times \mathbb D^{n_k} \subset \C^n$. This map is obviously $G$-equivariant, so passing to the quotients we obtain a homeomorphism $\hat \varphi: U/G \xrightarrow{\cong} SP^{n_1}(\mathbb D) \times \ldots \times SP^{n_k}(\mathbb D)\subset SP^{n_1}(\C) \times \ldots \times SP^{n_k}(\C)$.

Now consider the canonical homeomorphisms $h_n\colon SP^n(\C) \cong \C^n$, sending $[z_1, \ldots, z_n]$ to $(\sigma_1(z), \ldots, \sigma_n(z))$, where $\sigma_i$ is the $i$-th elementary symmetric polynomial in $z=(z_1, \ldots, z_n)$.

So, we obtain a map $$\Phi\colon V \xrightarrow[\cong]{\hat p^{-1}} U/G \xrightarrow[\cong]{\hat\varphi} \prod SP^{n_i}(\mathbb D) \subset \prod SP^{n_i}(\C) \xrightarrow[\cong]{\prod h_{n_i}} \prod \C^{n_i}=\C^n $$ This is a homeomorphism onto its image (which is open in $\C^n$); in other words, it is a (topological) chart around an arbitrary $x\in SP^n(M)$. Such charts obviously form an open cover of $SP^n(M)$ and are indexed by the following data: disjoint open sets $U_1,\ldots,U_k \subset M$, Klein charts $\varphi_i\colon U_i \to \mathbb D$, and multiplicities $n_1,\ldots, n_k\ge 1$.

Explicitly, consider any point $y=[y_1,\ldots, y_n]\in V$. By definition, exactly $n_i$ of the points $y_1, \ldots, y_n$ lie in $U_i$. There exists a unique monic complex polynomial $P_i(t)$ (of degree $n_i$) whose roots are exactly the $\varphi_i$-coordinates of these points. Then $\Phi$ sends $y$ to the coefficients of the polynomials $P_1, \ldots, P_k$.
\smallskip

2) \emph{Compatibility of charts}

a) Now we claim that the constructed charts are dianalytically compatible. In other words, let $(V, \Phi)$ and $(V', \Phi')$ be two such charts, constructed from $(U_i, \varphi_i, n_i)_{i=1}^k$ and $(U_j', \varphi'_j, n'_j)_{j=1}^s$ respectively. Then, if $W=V\cap V'\ne\varnothing$, we need to show that the composition $\Phi' \circ \Phi^{-1}: \Phi(W) \to \Phi'(W)$ is holomorphic or antiholomorphic (on each connected component).

Denote $W_{ij} = U_i \cap U'_j\subset M$. These are pairwise disjoint open subsets of $M$ (possibly empty).

For a point $y =[y_1,\ldots, y_n] \in W$ denote by $a_{ij}$ the number of points $y_i$ lying in $W_{ij}$. We obviously have $\sum_{i=1}^k a_{ij} = n'_j, \hspace{1ex}\sum_{j=1}^s a_{ij} = n_i, \hspace{1ex} \sum_{i,j} a_{ij} = n$.

For a fixed matrix $A = (a_{ij})$ of non-negative integers satisfying the above identities, consider an open subset $W_A = \{ y \in W \mid \text{the number of points in } W_{ij} \text{ is } a_{ij} \text{ for all } i,j \}$.

It is clear that $W=\bigsqcup_A W_A$. So, it is sufficient to check that $\Phi'\circ\Phi^{-1}$ is dianalytic on each $\Phi(W_A)$.
\smallskip

b) Consider nonempty $W_A$ and the corresponding $W_{ij}\subset M$ for which $a_{ij} > 0$. Notice that the data consisting of $W_{ij}$, $\varphi_i|_{W_{ij}}\colon W_{ij}\xrightarrow{\cong} \varphi_i(W_{ij})\subset \C$ and multiplicities $a_{ij}$ can be used to define the chart on $W_A$ by a procedure similar to the one above. Namely, we have $$\Psi\colon W_A\xrightarrow[\cong]{\hat p^{-1}} (\prod W_{ij}^{\times a_{ij}})/G_A\xrightarrow[\cong]{\hat\varphi}\prod SP^{a_{ij}}(\varphi_i(W_{ij}))\hookrightarrow \prod SP^{a_{ij}}(\C )\cong \prod \C^{a_{ij}}=\C^n$$ This is a homeomorphism onto its image.

Consider the composition $F = \Phi \circ \Psi^{-1}: \Psi(W_A) \to \Phi(W_A)$. By definition, for fixed $i$, the chart $\Psi$ gives, for each $j$, a monic polynomial $P_{ij}(t)$ of degree $a_{ij}$ whose roots are the $\varphi_i$-coordinates of the points in $W_{ij}$. Similarly, the chart $\Phi$ requires a single monic polynomial $P_i(t)$ of degree $n_i$ whose roots are the $\varphi_i$-coordinates of all points in $U_i$. Since the $W_{ij}$ (for fixed $i$) are disjoint and their union is $U_i$ (on the domain where the points lie exactly in these sets), we obtain that $P_i(t) = \prod_{j} P_{ij}(t)$. Since polynomial multiplication is a polynomial map on the coefficients, we have that $F$ is a polynomial map, hence holomorphic.

Now consider the composition $G = \Phi' \circ \Psi^{-1}: \Psi(W_A)\to \Phi'(W_A)$. It is given as follows. For each $(i,j)$, the chart $\Psi$ gives the polynomial $P_{ij}(t)$ with roots equal to the $\varphi_i$-coordinates of the points in $W_{ij}$. On the other hand, if we have the polynomial $Q_{ij}(t)$ whose roots are the $\varphi'_j$-coordinates of the same points, then the chart $\Phi'$ requires the polynomial $Q_j = \prod Q_{ij}$. So we need to check that the coefficients of $Q_{ij}$ are dianalytic functions of the coefficients of $P_{ij}$.

Let $a = a_{ij}$. Denote the roots of $P_{ij}$ by $z_1,\ldots,z_a \in \varphi_i(W_{ij})$. Then the roots of $Q_{ij}$ are $\varphi_{ij}(z_1),\ldots,\varphi_{ij}(z_a)$, where $\varphi_{ij} = \varphi'_j \circ \varphi_i^{-1}: \varphi_i(W_{ij}) \to \varphi'_j(W_{ij})$ is a transition function on $M$, so it is either a biholomorphism or an anti-biholomorphism between open subsets of $\mathbb{D}$.

The coefficients of $Q_{ij}$ are the elementary symmetric polynomials $\sigma_k(\varphi_{ij}(z_1),\ldots,\varphi_{ij}(z_a))$. This is a symmetric function in $(z_1,\ldots,z_a)$ on the $\Sigma_a$-invariant open set $\varphi_i(W_{ij})^{\times a} \subset \mathbb{C}^a$.

If $\varphi_{ij}$ is holomorphic, then $\sigma_k(\varphi_{ij}(z_1),\ldots,\varphi_{ij}(z_a))$ is also holomorphic. Hence, by the Holomorphic Symmetric Function Theorem, we have $\sigma_k(\varphi_{ij}(z_1),\ldots,\varphi_{ij}(z_a))=f(\sigma_1(z), \ldots, \sigma_a(z))$ for a (unique) holomorphic $f$. But $\sigma_k(z)$ are exactly the coefficients of $P_{ij}$, so we obtain that the coefficients of $Q_{ij}$ depend holomorphically on the coefficients of $P_{ij}$, as claimed. For antiholomorphic $\varphi_{ij}$ we have that $\overline{\sigma_k(\varphi_{ij}(z_1),\ldots,\varphi_{ij}(z_a))}=\sigma_k(\overline{\varphi_{ij}}(z_1), \ldots, \overline{\varphi_{ij}}(z_a))$ is holomorphic and symmetric, so it can also be written as $f(\sigma_1(z), \ldots, \sigma_a(z))$, and in this case the coefficients of $Q_{ij}$ depend antiholomorphically on the coefficients of $P_{ij}$.

c) Thus, we have proved that $F = \Phi \circ \Psi^{-1}$ and $G = \Phi' \circ \Psi^{-1} $ are dianalytic. Now we need to check that $\Phi' \circ \Phi^{-1} = G \circ F^{-1}$ is also dianalytic (on $\Phi(W_A)$).

Recall that $F$ is just a polynomial map given by multiplying polynomials. We claim that it is a biholomorphism onto its image.

Indeed, $F$ is a product over $i$ of maps $$
F_i: \prod_j SP^{a_{ij}}(\varphi_i(W_{ij})) \to SP^{n_i}(\varphi_i(U_i)), \quad (P_{ij})_j \mapsto \prod_j P_{ij}.
$$
Here we identify $SP^a(W)$ with the space of monic polynomials of degree $a$ with roots in $W$.

The map $F_i$ is injective by the unique factorization of polynomials over $\C$ (since the roots of $P_{ij}$ lie in the disjoint sets $W_{ij}$).

The image of $F_i$ is the set of polynomials of degree $n_i$ that have exactly $a_{ij}$ roots in each $W_{ij}$ (with multiplicities). We need to construct a holomorphic inverse. In other words, we need to check that the coefficients of the factors in the unique factorization depend holomorphically on the coefficients of the product. But this is standard, because for any polynomial $P$ and any subset $S$ of its roots, if we choose a closed contour $\gamma$ that encloses all roots in $S$ (and no other roots), then by the residue theorem, the contour integral $s_r = \frac{1}{2\pi i} \oint_\gamma z^r \frac{P'(z)}{P(z)} \, dz$ gives an $r$-th power sum of the roots in $S$ (depending complex-analytically on the coefficients of $P$). Then, by Newton's identities, we can reconstruct the elementary symmetric polynomials in these roots, thus obtaining the corresponding factor of $P$.

So, $F^{-1}$ is holomorphic. Hence, the atlas we constructed on $SP^n(M)$ is dianalytic, as claimed. The functoriality with respect to Klein maps can easily be checked in a similar way (if we consider the local coordinate expression of a dianalytic map instead of the transition functions $\varphi_{ij}$).

If all transition functions on $M$ are holomorphic, then the same is true for $SP^n(M)$.
\end{proof}

\begin{remark}
    One can check that the constructed dianalytic (or complex) structure is universal in the following sense: for any dianalytic (complex) $N^{2n}$, the projection $p\colon M^{\times n} \to SP^n(M)$ induces a bijection between dianalytic (holomorphic) maps $SP^n(M)\to N$ and $\Sigma_n$-invariant dianalytic (holomorphic) maps $M^{\times n} \to N$.
\end{remark}

\begin{remark}
    Notice that the corresponding smooth structure is not functorial with respect to arbitrary smooth maps. For example, one can check that for the diffeomorphism $\R^2 \to \R^2$, $(x,y) \mapsto (x, 2y)$, the induced map $\R^4\cong SP^2({\R^2}) \to SP^2(\R^2) \cong \R^4$ is non-smooth at the origin.

    In particular, dianalytic coordinates on $M$ are essential and the smooth structure on $M$ alone is not enough to construct a smooth structure on $SP^n(M)$.
\end{remark}

\begin{remark}\label{BGZrem}
    In our proof of the (non-orientable version of the) Blagojevi\'{c}--Gruji\'{c}--Zivaljevi\'{c} conjecture (Theorem \ref{BGZ}) we considered the inclusion $N_{g,s}\hookrightarrow N_g$ and then used that the induced open embedding $SP^n(N_{g,s})\hookrightarrow SP^n(N_g)$ is smooth. This is true for the induced dianalytic structure on $N_{g,s}\subset N_g$. Since in Theorem \ref{BGZ} we are essentially interested in the topological type of the symmetric product, it does not matter which dianalytic (and hence smooth) structure we consider on $SP^n(N_{g,s})$.
\end{remark}

\begin{remark}
    One can check that the constructed smooth structure can also be characterized as follows: it is the unique smooth structure on symmetric products such that

1) it is dianalytically functorial;

    2) the projection $p\colon M^{\times n} \to SP^n(M)$ is smooth;

    3) on $SP^n(\C)$ it coincides with the one constructed from the elementary symmetric polynomial identification $SP^n(\C) \cong \C^n$.
\end{remark}

\begin{remark}

Consider the singular locus $\Delta=\{(x_1, \ldots, x_n) \in M^{\times n} \mid \exists\, i \ne j \mbox{ s.t. } x_i=x_j \}$. The action $\Sigma_n \curvearrowright M^{\times n}$ is free outside $\Delta$, so we obtain a canonical smooth structure on an open dense subset $U_{reg}=(M^{\times n}-\Delta)/\Sigma_n \subset SP^n(M)$. One can check that this canonical smooth structure on $U_{reg}$ coincides with the restriction of the constructed smooth structure on $SP^n(M)$.

In other words, $p\colon M^{\times n} \to SP^n(M)$ is a smooth branched covering map ramified over $p(\Delta)=SP^n(M)-U_{reg}$.
\end{remark}

\begin{remark}
        By the classical Chevalley--Shephard--Todd theorem \cite{st}, \cite{ch}, for any finite subgroup $G\subset GL(n, \C)$ the natural quotient analytic structure on $\C^n/G$ is smooth if and only if $G$ is generated by complex pseudoreflections. In the case of subgroups of $\Sigma_n$ this is equivalent to being generated by transpositions.
    \end{remark}

\begin{example}
    One can check that the natural bijection $SP^n(\C P^1)\cong \C P^n$, sending $[(z_1:w_1),\ldots, (z_n:w_n)]$ to the coefficients of the polynomial $\prod_{i=1}^n(z_i t - w_i)$, provides a biholomorphism between $SP^n(\C P^1)$ with the constructed complex structure and $\C P^n$ with its standard complex structure.

    In particular, we see that in Theorem \ref{th1} a) we have the standard complex projective spaces.

\end{example}

\begin{example}
    Consider the free antipodal antiholomorphic involution $\tau$ on $\C P^1$. Then it obviously extends to $SP^{2n}(\C P^1)=\C P^{2n}$, and as we mentioned above the fixed-point set $SP^{2n}(\C P^1)^{\tau}$ can be identified with $SP^n(\C P^1/\tau)=SP^n(\R P^2)$.

    On the other hand, for this explicit involution one can directly check that $SP^{2n}(\C P^1)^{\tau}=(\C P^{2n})^\tau$ can be identified with $\R P^{2n}$ (with the standard real-analytic structure). From this, one can obtain that $SP^n(\R P^2)$ with the constructed above smooth structure coincides with $\R P^{2n}$ with the standard smooth structure.

    This is exactly the smooth structure considered by Dupont and Lusztig in \cite{dulu}. So, we see that in Theorem \ref{th1} b) we have the standard real projective spaces.
\end{example}

\begin{remark}\label{taufix}
 In general, let a nonorientable surface $M$ be represented as a quotient $\widetilde M/\tau$ of a complex curve $\widetilde M$ by a free antiholomorphic involution $\tau$. Then for even $n$, we have $SP^n(\widetilde M)^\tau=SP^{n/2}(M)$. In this case, the standard real-analytic structure on the fixed-point set of an antiholomorphic involution on a complex manifold coincides with the one induced by the dianalytic structure on $SP^{n/2}(M)$ constructed above.
 
 On the other hand, for odd $n$, the extended involution $\tau$ is free on $SP^n(\widetilde M)$. Notice that the quotient $SP^n(\widetilde M)/\tau$ is not the same as $SP^n(M)$ (not even homotopy equivalent).
\end{remark}

While in general the constructed dianalytic structure on $SP^n(M)$ depends on the dianalytic structure on $M$, it turns out that the corresponding smooth structure does not.

\begin{proposition}\label{smun}
    Let $M$ be a closed surface. Then the smooth structure on $SP^n(M)$ constructed above does not depend (up to a diffeomorphism isotopic to the identity) on the choice of the Klein structure on $M$.
\end{proposition}

\begin{remark}
    Notice that, in general, there exist exotic smooth structures on $SP^n(M)$. For example, there are exactly eight smooth structures on $\C P^5=SP^5(\C P^1)$, and for $6 \le m \le 8$ there are exactly four smooth structures on $\C P^m=SP^m(\C P^1)$ (up to diffeomorphism), see \cite{kas26}.
\end{remark}

Before proving this, we need some preliminary observations. First of all, notice that, for any continuous map $E\to B$, we can consider the relative Cartesian $n$-fold product $\underbrace{E \times_B \times \cdots \times_B E}_{n \mbox{ {\scriptsize times}}}\to B$. Since there is a natural permutation action on it, we can also naturally define the \emph{relative symmetric product} $$SP^n(E/B):=\Bigr(\underbrace{E \times_B \times \cdots \times_B E}_{n \mbox{ {\scriptsize times}}}\Bigl)/\Sigma_n \to B$$ Notice that, if the map $f$ is proper (i.e., preimages of compact substes are compact), then the map $SP^n(E/B) \to B$ is also proper.

Secondly, by a \emph{smooth family of Klein surfaces} we mean a smooth submersion $f\colon E \to B$ with connected $B$ satisfying the following properties:

1) there is an open cover of $E$ by charts $W_\alpha$ of $E$ of the form $\Phi_\alpha\colon W_\alpha \xrightarrow{\cong} D_\alpha \times U_\alpha$, where $D_\alpha \subset \C$, $U_\alpha \subset B$, and $pr_2 \circ \Phi_\alpha = f|_{W_\alpha}$; in particular, we obtain an atlas $\Phi_{\alpha, b} \colon W_\alpha \cap E_b \xrightarrow{\cong} D_\alpha$ on each fiber $E_b=f^{-1}(b)$, $b\in B$;

2) the atlases above must be dianalytic, that is, the compositions $\Phi_{\alpha} \circ \Phi_{\beta}^{-1} \colon \Phi_{\beta}(W_\alpha \cap W_\beta) \to \Phi_{\alpha}(W_\alpha \cap W_\beta)$ must have the form $(z, b) \mapsto (\varphi_{\alpha \beta}(z), b)$, where $\varphi_{\alpha\beta}$ is either holomorphic or antiholomorphic.

Now, if $f\colon E \to B$ is a smooth family of Klein surfaces, then one can directly check that the above construction of a smooth structure on symmetric products works as well in the relative situation and gives a smooth-submersion structure on $SP^n(E/B) \to B$ with fibers $SP^n(E_b)$ (constructed from the corresponding dianalytic structure on $E_b$).

If a family $f\colon E \to B$ is proper, then the submersion $SP^n(E/B) \to B$ is also proper. But by Ehresmann's theorem, any proper submersion is a smooth fiber bundle. In particular, all fibers $SP^n(E_b)$ are diffeomorphic, and a diffeomorphism can be chosen to be isotopic to the identity.

The upshot is that, if we can include two dianalytic structures on $M$ in a proper smooth family, then the corresponding smooth structures on symmetric products are diffeomorphic (by a diffeomorphism isotopic to identity).

\begin{proof}[Proof of Proposition \ref{smun}]

\emph{1)} Let $\mathcal A_1$ and $\mathcal A_2$ be two dianalytic structures on a closed surface $M$. They define two smooth structures on $M$: $\mathcal S_1$ and $\mathcal S_2$. They need not coincide, but they are diffeomorphic (as any two smooth structures on a surface are). Consider any diffeomorphism $f\colon (M, \mathcal S_1) \xrightarrow{\cong} (M, \mathcal S_2)$. We may also assume that $f$ is isotopic to identity. We can consider a pull-back $f^*(\mathcal A_2)$ of the dianalytic structure $\mathcal{A}_2$, so by definition $f$ becomes a dianalytic isomorphism $f\colon (M, f^*(\mathcal A_2)) \xrightarrow{\cong} (M, \mathcal A_2)$.

But since our construction of a smooth structure on symmetric products is functorial with respect to dianalytic maps, we obtain that $SP^n(f)\colon SP^n(M, f^*(\mathcal A_2)) \to SP^n(M, \mathcal A_2)$ is a diffeomorphism (isotopic to identity, since $f$ is).

On the other hand, now $f^*(\mathcal A_2)$ and $\mathcal A_1$ define the same smooth structure $\mathcal S_1$ on $M$. So, now we need to check that two dianalytic structures defining the same smooth structure on $M$ define diffeomorphicsmooth structures on symmetric products (via a diffeomorphism isotopic to identity).

\emph{2) Oriented case.} Let $M$ be oriented. Let us prove that the smooth structure on symmetric products does not depend on the choice of a complex structure on $M$ compatible with the fixed orientation and smooth structure on $M$.

Recall that for any (almost) complex structure $J$ (compatible with the orientation) there exists a unique Riemannian metric $g$ such that $g(JX, JY)=g(X,Y)$. On the other hand, for a given $g$ there exist only two such complex structures: $J$ and $-J$, and only $J$ is compatible with the orientation. Clearly, the bijection $g \leftrightarrow J$ is smooth.

Let $J_0$ and $J_1$ be two almost complex structures on an oriented $M$. Then we can connect the corresponding metrics by a linear path $g_t = (1-t) g_0 + t g_1$. So, we obtain a smooth path $J_t$ of almost complex structures between $J_0$ and $J_1$. To obtain a smooth family of complex Riemann surfaces (of the form $M\times[0,1]\to [0,1]$) we need only be able to smoothly integrate these almost complex structures on $M$. But it is well known that the corresponding local complex coordinate $z_t$ is uniquely determined by its real part ${\rm Re}\, z_t=u_t$, which is a harmonic function with respect to the Laplace--Beltrami operator of $g_t$, so it depends smoothly on $t$ (since the PDE under consideration is elliptic). So we obtain a smooth family of complex Riemann surfaces of the form $E=M\times [0,1] \to [0,1]$.

So, by the relative symmetric product construction and Ehresmann's theorem we obtain a fiber bundle $SP^n(E/[0,1])\to [0,1]$. Hence $SP^n(M, J_0)$ is diffeomorphic to $SP^n(M, J_1)$. Moreover, Ehresmann's theorem gives a diffeomorphism isotopic to the identity (here we assume that all the complex structures are compatible with the orientation).

\emph{3) General case.} Consider the orientation double cover $\widetilde M \to M$ (trivial, if $M$ is orientable) with the orientation-reversing deck involution $\tau$. Any dianalytic structure on $M$ can be lifted uniquely to a complex structure on $\widetilde M$ compatible with the orientation and such that $\tau$ is antiholomorphic with respect to it.

So let $J_0$ and $J_1$ be two almost complex structures on $\widetilde M$ corresponding to dianalytic structures $\mathcal A_0$ and $\mathcal A_1$ on $M$. Since $\tau$ is antiholomorphic, we have $\tau^*(J_i)=-J_i$. But this means that for the corresponding metrics we have $\tau^*g_i=g_i$, i.e., $g_i$ are $\tau$-invariant.

Then $g_t=(1-t)g_0+tg_1$ is also $\tau$-invariant. Hence, the corresponding $J_t$ are anti-$\tau$-invariant (i.e., $\tau^*J_t=-J_t$). So, $\tau$ is antiholomorphic with respect to all the corresponding complex structures.

Again, we have a smooth family of complex Riemann surfaces $\widetilde M\times [0,1]\to [0,1]$. Taking the quotient $(\widetilde M\times [0,1])/(\tau\times \rm id)$ we obtain the required smooth family of Klein surfaces of the form $M\times [0,1]\to [0,1]$. Then again the relative symmetric product and Ehresmann's theorem give the required diffeomorphism isotopic to the identity.
\end{proof}

\begin{remark}
    While the complex structure on $SP^n(M)$ may depend on the choice of a complex structure on $M$, we see that the corresponding almost complex structures are homotopic (since they are homotopic on $M$). In particular, the stable complex structure on $SP^n(M)$ is also canonical (this holds for any surface). The resulting complex bordism class is given by Theorem \ref{th1} b).
\end{remark}

\begin{remark}
   Of course, for any surface $M$ the constructed family $SP^n(M\times[0,1]/[0,1])$ is topologically trivial, i.e., homeomorphic to $SP^n(M)\times[0,1]\to [0,1]$ (notice that to construct it we do not need the compactness of $M$). However, properness and Ehresmann's theorem are needed to deduce a \emph{smooth} trivialization (since a priori the smooth structure on $SP^n(M\times[0,1]/[0,1])$ is not the product one).
\end{remark}
   
One can ask whether Proposition \ref{smun} holds for noncompact surfaces.

Since any noncompact surface can be exhausted by surfaces of finite type, it is essentially sufficient to consider the latter.

A surface of finite type is always diffeomorphic to a compact surface $\overline M$ with finitely many punctures. But dianalytically this is not true: it can have dianalytically nontrivial “boundaries”.

If we consider dianalytic structures on a surface $M$ that are induced from dianalytic structures on $\overline M$ (in other words, structures that can be extended over each puncture), then Proposition \ref{smun} holds. Indeed, let $\mathcal A_0$ and $\mathcal A_1$ be two dianalytic structures on $M$ induced from dianalytic structures $\overline{\mathcal A_i}$ on $\overline{M}$. Since $\overline{\mathcal A_0}$ and $\overline{\mathcal A_1}$ are equivalent locally near the marked points, there exists a diffeomorphism $f\colon \overline M\to \overline M$, isotopic to the identity and preserving the marked points throughout the isotopy, such that $f^*(\overline {\mathcal A_1})$ coincides with $\overline {\mathcal A_0}$ near the marked points. Since
$f\colon (\overline M, f^*(\overline {\mathcal A_1}))\to (\overline M, \overline {\mathcal A_1})$
is a dianalytic isomorphism, it induces a diffeomorphism (isotopic to identity) on symmetric products. Since it preserves the punctures, it also induces a diffeomorphism on the symmetric products of $M$. Thus we may assume that $\overline {\mathcal A_1}$ and $\overline {\mathcal A_0}$ coincide near the punctures.

Then the path between these two dianalytic structures on $\overline{M}$ can be chosen to be constant near the punctures. Hence the resulting family of Klein surfaces over the interval can be chosen to be trivial near the punctures. But then it is also trivial after passing to symmetric products. Therefore Ehresmann’s diffeomorphism $SP^n(\overline M, \overline{\mathcal{A}_0})\cong SP^n(\overline M, \overline{\mathcal{A}_1})$ can be chosen to be the identity near the marked points, and hence restricts to a diffeomorphism $SP^n(M, \mathcal{A}_0)\cong SP^n(M, \mathcal A_1)$.

However, this argument does not work for arbitrary surfaces, since there are dianalytic structures that cannot be extended to closed surfaces (“has non-extendable ends”), so Ehresmann’s theorem cannot be applied here.

On the other hand, using much more nontrivial machinery, one can check that the following actually holds.

\begin{proposition}
    Let $N$ be a topological manifold (without boundary) of dimension $\dim N\ge 5$. Assume that we have a smooth structure on $N\times [0,1]$ such that the projection $N\times [0,1]\to [0,1]$ is a smooth submersion. Then the smooth structures on $N\times 0$ and $N\times 1$ are diffeomorphic via a diffeomorphism isotopic to the identity.

In particular, this holds for the relative symmetric product constructed above from a smooth family of Klein structures on an arbitrary surface. Hence we obtain that the smooth structure on $SP^n(M)$ indeed does not depend on the Klein structure on $M$, if $n>2$.
\end{proposition}

\begin{proof}
    This is a direct consequence of the smoothing theory of Kirby and Siebenmann \cite{ks}. The central classification theorem of Essay V is the following.
    
\begin{fact}[{\cite[Essay V, Theorem 2.3]{ks}}]
    Let $N$ be a topological manifold (without boundary) of dimension $n=\dim N \ne 4$. Then there exists a homotopy equivalence between the spaces ${\rm DIFF}(N)$ and ${\rm Lift}(TN, BO(n))$.
\end{fact}

Here ${\rm DIFF}(N)$ and ${\rm Lift}(TN, BO(n))$ are (geometric realizations of) the following semi-simplicial sets. For ${\rm DIFF}(N)$, a $k$-simplex is a smooth structure on $N\times \Delta^k$ such that the projection $N\times \Delta^k\to \Delta^k$ is a smooth submersion. In particular, in our case we have, by definition, a path (an edge) $[0,1]\to {\rm DIFF}(N)$.

Let $f\colon N\to BTop(n)$ be a (fixed) classifying map for a topological tangent (micro)bundle $TN$. Then a $k$-simplex of ${\rm Lift}(TN, BO(n))$ is a map $\sigma\colon \Delta^k \times N \to BO(n)$ such that the following diagram commutes:
\[
\begin{tikzcd}
    \Delta^k\times N \arrow[d] \arrow[r, "\sigma"] & BO(n)\arrow[d]\\
    N \arrow[r, "f"] & BTop(n)
\end{tikzcd}
\]
(where we assume that $BO(n)\to BTop(n)$ is a fibration).

Now, the central classification theorem of Essay IV is the following.
\begin{fact}[{\cite[Essay IV, Theorem 10.1]{ks}}]
    Let $N$ be a topological manifold (without boundary) of dimension $n=\dim N\ge 5$. Then there exists a bijection between the sets $\mathcal S_{\rm DIFF}(N)$ and $\pi_0({\rm Lift}(TN, BO(n)))$.
\end{fact}

Here $\mathcal S_{\rm DIFF}(N)$ is the set of smooth structures on $N$ modulo concordance. As a result, in our situation we obtain that smooth structures on $N\times 0$ and $N\times 1$ correspond to the same path component of ${\rm Lift}(TN, BO(n))$, and hence they are concordant. Finally, by \cite[Essay I, Theorem 4.1]{ks} (if $\dim N\ge 5$) concordant smooth structures on $N$ are isotopic.
\end{proof}

\begin{remark}
    Since we have continuous families of exotic smooth structures on open $4$-dimensional manifolds, the question of whether the smooth structure on $SP^2(M)$ is independent of the Klein structure on $M$ seems a little bit more mysterious (for noncompact $M$). It would be interesting to understand whether this is true.
\end{remark}


\begin{thebibliography}{99}

\bibitem[AG]{ag} N. L. Alling, N. Greenleaf, \emph{Foundations of the theory of Klein surfaces}. Lecture Notes in Mathematics, Vol. 219., Springer-Verlag, 1971.

\bibitem[Arn96]{arn} V. Arnold, \emph{Topological content of the Maxwell theorem on multipole representation of spherical functions}. Topol. Methods Nonlinear Anal., 7:2 (1996), 205--217.

\bibitem[BGZ03]{bgz1} P. Blagojevi\'{c}, V. Gruji\'{c}, R. Zivaljevi\'{c}, \emph{Symmetric products of surfaces and the cycle index}, Israel J. Math. 138 (2003), 61-72.

\bibitem[BGZ05]{bgz2} P. Blagojevi\'{c}, V. Gruji\'{c}, R. Zivaljevi\'{c}, \emph{Arrangements of symmetric products of spaces}, Topology and its Applications 148 (2005), 213-232.

\bibitem[BGZ]{bgz} P. Blagojevi\'c, V. Gruji\'c, R. Zivaljevi\'c, \emph{Symmetric products of surfaces; a unifying theme for topology and physics}, arXiv:math/0408417.

\bibitem[BR18]{br} Y. Boote, N. Ray, \emph{On the symmetric squares of complex and quaternionic projective space}. Glasgow Mathematical Journal, 60(3), 633–648, (2018).

\bibitem[Buch06]{buch06} V. M. Buchstaber, \emph{$n$-valued groups: theory and applications}, Mosc. Math. J., 6:1 (2006), 57–84.

\bibitem[Buch12]{buch12} V. M. Buchstaber, \emph{Complex cobordism and formal groups}, Russian Math. Surveys, 67:5 (2012), 891–950.

\bibitem[BP]{bp} V. M. Buchstaber, T. E. Panov, \emph{Toric Topology}, Mathematical Surveys and Monographs, vol. 204, American Mathematical Society, Providence, RI, 2015.

\bibitem[Ch55]{ch} C. Chevalley, \emph{Invariants of finite groups generated by reflections}, Amer. J. Math., 77 (4), 778–782, (1955).

\bibitem[CF64]{cf} P. E. Conner, E. E. Floyd, \emph{Differentiable Periodic Maps}, Academic Press Inc., New York, 1964.

\bibitem[Dol58]{dold} A. Dold, \emph{Homology of Symmetric Products and Other Functors of Complexes}. The Annals of Mathematics, Second Series, Vol. 68, No. 1 (Jul., 1958), pp. 54-80.

\bibitem[DP61]{dp} A. Dold., D. Puppe, \emph{Homologie nicht-additiver Funktoren. Anwendungen.} Annales de l'Institut Fourier, 11, (1961), 201--312.

\bibitem[DT58]{dt} A. Dold, R. Thom, \emph{Quasifaserungen und unendliche symmetrische Produkte}. Annals of Mathematics, 67(2), (1958), 239–281.

    \bibitem[DL71]{dulu} J. L. Dupont, G. Lusztig, \emph{On manifolds satisfying $w^2_1 = 0$}. Topology, 10:2 (1971), 81--92.

    \bibitem[EGL01]{egl} G. Ellingsrud, L. G\"ottsche, M. Lehn, \emph{On the cobordism class of the Hilbert scheme of a surface}. J. Algebraic Geom. 10 (2001), no. 1, 81--100.

    \bibitem[Ep66]{eps} D. B. A. Epstein, \emph{Curves on $2$-manifolds and isotopies}, Acta Math. 115, 83–107 (1966).

    \bibitem[Fo68]{fog} J. Fogarty, \emph{Algebraic Families on an Algebraic Surface}, American Journal of Mathematics , Apr., 1968, Vol. 90, No. 2 (Apr., 1968), pp. 511--521.

\bibitem[Gu11]{gug1} D. V. Gugnin, \emph{Topological applications of graded Frobenius $n$-homomorphisms}, Trans. Moscow Math. Soc., 2011, 97–142.

\bibitem[Gu12]{gug2} D. V. Gugnin, \emph{Topological applications of graded Frobenius $n$-homomorphisms, II}, Trans. Moscow Math. Soc., 2012, 167–182.

\bibitem[Gu22]{gug3} D. V. Gugnin, \emph{Any Suspension and Any Homology Sphere Are
2H-Spaces}, Proc. Steklov Inst. Math., 318 (2022), 45–58.

    \bibitem[Gu24]{gug4} D. V. Gugnin, \emph{The Integral Cohomology Ring of Symmetric Products of CW Complexes and Topology of Symmetric Products of Riemann Surfaces}, Proc. Steklov Inst. Math., 326 (2024), 133--156.

\bibitem[GGK]{ggk} V. Guillemin, V. Ginzburg, Y. Karshon, \emph{Moment Maps, Cobordisms, and Hamiltonian Group Actions}, Mathematical Surveys and Monographs 98 (2002).

\bibitem[Ha]{ha26} A. Hatcher, \emph{The Kirby torus trick for surfaces},
    Enseign. Math. (2) 72 (2026), 161–174.

\bibitem[GR]{gr} R. C. Gunning, H. Rossi, \emph{Analytic Functions of Several Complex Variables}. Prentice-Hall, 1965.


\bibitem[Ja25]{j25} E. Jauhari, \emph{LS-category and sequential topological complexity of symmetric products}, J. Appl. Comput. Topol. 9, No. 4, Paper No. 25, 31 p. (2025).

\bibitem[KT10]{kt} S. Kallel, W. Taamallah, \emph{Fundamental Group and Euler Characteristic of Permutation Products and Fat Diagonals}, arXiv.org/abs/1010.1507.

    \bibitem[Kas26]{kas26} R. Kasilingam, \emph{Diffeomorphism Classification of Smooth Structures and Tangential Homotopy Types of $\C P^m$ for $5\le m\le 8$}, arXiv:2604.27521.

\bibitem[Ke01]{ke01} S. Kera, \emph{On the permutation products of manifolds}, Beiträge zur Algebra und Geometrie (Contr. to Alg. and Geom.), Vol. 42, 2, p. 547–555, 2001.

\bibitem[KS77]{ks} R. C. Kirby, L. C. Siebenmann, \emph{Foundational Essays on Topological Manifolds, Smoothings, and Triangulations}, Annals of Mathematics Studies, No. 88,
Princeton University Press, Princeton, N.J., 1977.

    \bibitem[La19]{la} C. Lange, \emph{When is the underlying space of an orbifold a manifold?}, Trans. of the Amer. Math. Soc., vol. 371, no. 9, pp. 6633–6652 (2019).

    \bibitem[Mcd62]{mcd} I. G. Macdonald, \emph{Symmetric product of an algebraic curve}. Topology, 1:4 (1962), 319--343.

    \bibitem[Mat61a]{mat1} A. Mattuck, Picard bundles. Illinois J. Math., 5 (1961), 550--564.

\bibitem[Mat61b]{mat2} A. Mattuck, \emph{Symmetric Products and Jacobians.}, American Journal of Mathematics, vol. 83, no. 1, 1961, 189--206.

 \bibitem[Mat62]{mat3} A. Mattuck, \emph{On Symmetric Products of Curves}, Proc. Amer. Math. Soc. 13 (1962), 82--87.

\bibitem[Milg69]{milg} R. J. Milgram, {The homology of symmetric products.} Transactions of the American Mathematical Society, 138, (1969), 251–265.

 \bibitem[Miln60]{mil} J. Milnor, \emph{On the cobordism ring $\varOmega_*$ and a complex analogue, Part I}, American Journal of Mathematics, 82 (3): 505–-521, (1960).

 \bibitem[Nak57]{nak} M. Nakaoka, \emph{Cohomology of symmetric products}, J. Inst. Polytechn., Osaka City Univ., Ser. A 8, 121-145 (1957).

 \bibitem[Nov62]{nov} S. P. Novikov, \emph{Homotopy properties of Thom complexes}, Mat. Sb., New Series, 57: 407–-442, (1962) (in Russian).

 \bibitem[Nov67]{nov2}  S. P. Novikov, \emph{The methods of algebraic topology from the viewpoint of cobordism theory}, Math. USSR-Izv., 1:4 (1967), 827–-913.

\bibitem[Ong03]{ong} B. W. Ong \emph{The homotopy type of the symmetric products of bouquets of circles}, Int. J. Math. 2003. V. 14, N 5. P. 489–497.

\bibitem[Ra25]{rado} T. Rad\'o, \emph{\"Uber den Begriff der Riemannschen Fl\"ache}, Acta Sci. Math. Szeged. 2 (1925), 101–121.

\bibitem[Rav]{rav}
D. C. Ravenel, \emph{Complex Cobordism and Stable Homotopy Groups of Spheres},
Pure and Applied Mathematics,~121, Academic Press Inc., Orlando, FL, 1986.

   \bibitem[Sch63]{sch} R. L. E. Schwarzenberger, \emph{Jacobians and symmetric products}, Illinois J. Math. 7 (2), 257--268, (1963).

   \bibitem[ST54]{st} G. S. Shephard, J. A. Todd, \emph{Finite unitary reflection groups}, Can. J. Math., 6, 274–304, (1954).

   \bibitem[Sto]{sto} R. E. Stong, \emph{Notes on cobordism theory}, Princeton University Press and the University of Tokyo Press, Princeton, New Jersey, 1968.

\bibitem[Th52]{thom} R. Thom, \emph{Espaces fibr\'es en sph\`eres et carr\'es de Steenrod}, Annales scientifiques de l'École Normale Supérieure, Serie 3, Volume 69 (1952), pp. 109--182.

\bibitem[Th54]{thom2} R. Thom, \emph{Quelques propriétés globales des variétés différentiables}, Comment. Math. Helv. 28, (1954), pp. 17--86

\bibitem[Wa80]{wa} C. H. Wagner, \emph{Symmetric, cyclic, and permutation products of manifolds}. Dissertationes Mathematicae, Vol. 182, pp. 1--52, (1980).

\end{thebibliography}
\end{document}